\documentclass[reqno,fleqn]{amsart}
\usepackage{acronym}
\usepackage{amsthm}
\usepackage[usenames,dvipsnames]{pstricks}
\usepackage{epsfig}
\usepackage{pst-grad} 
\usepackage{pst-plot} 
\usepackage[english]{babel}
\usepackage{graphicx,subfigure}
\usepackage{tabularx}
\usepackage{amsmath,amssymb}
\usepackage{algorithm,algpseudocode} 

\renewcommand{\P}{\mathcal{P}} 

\newcommand{\half}{{\textstyle{1\over2}}}

\makeatletter
\newcommand{\sech}{\mathop{\operator@font sech}}
\newcommand{\sign}{\mathop{\operator@font sign}}
\makeatother

\newtheorem{remark}{Remark}[section]

\acrodef{Serre}{Serre}
\acrodef{SGN}{Serre-Green-Naghdi}
\acrodef{cB}{`classical' Boussinesq}
\acrodef{FEM}{Finite Element Method}
\acrodef{IBVP}{initial-boundary value problem}
\acrodef{RK}{Runge-Kutta}
\acrodef{ODEs}{ordinary differential equations}
\acrodef{BT}{Boussinesq type}

\begin{document}

\title[A Modified Galerkin method for the SGN system]{A Modified Galerkin / Finite Element Method for the numerical solution of the Serre-Green-Naghdi system}

\author{Dimitrios Mitsotakis}
\address{Victoria University of Wellington, School of Mathematics and Statistics, PO Box 600, Wellington 6140, New Zealand}
\email{dimitrios.mitsotakis@vuw.ac.nz}
\urladdr{http://dmitsot.googlepages.com/}

\author{Costas Synolakis}
\address{University of Southern California, Viterbi School of Engineering, 3620 S. Vermont Avenue, 
Los Angeles, CA 90089, USA}
\email{costas@usc.edu}

\author{Mark McGuinness}
\address{Victoria University of Wellington, School of Mathematics and Statistics, PO Box 600, Wellington 6140, New Zealand}
\email{Mark.McGuinness@vuw.ac.nz}



\begin{abstract}
A new modified Galerkin / Finite Element Method  is proposed for the numerical solution of the  fully nonlinear shallow water wave equations. The new numerical method allows the use of low-order Lagrange finite element spaces, despite the fact that the system contains third order spatial  partial derivatives for the depth averaged velocity of the fluid. After studying the efficacy and the conservation properties of the new numerical method, we proceed with the validation of the new numerical model and  boundary conditions by comparing the numerical solutions with laboratory experiments and with available theoretical asymptotic results.
\end{abstract}

\maketitle


\section{Introduction}\label{introduction}

The motion of an ideal (inviscid, irrotational) fluid bounded above by a free surface and below by an impermeable bottom is governed by the full Euler equations of water wave theory, \cite{Whitham1999}. Because of the complexity of the Euler equations, a number of simplified models describing inviscid fluid flow have been derived such as various \acf{BT} models. The \acf{SGN} system can be considered to be a \acs{BT} model that approximates the Euler equations, and models one-dimensional, two-way propagation of long waves, without any restrictive conditions on the wave height. 
The \acs{SGN} system is a fully non-linear system of the form,
\begin{subequations} \label{eq:SGN}
\begin{align}
& h_t+(hu)_x =0\ , \label{eq:SGNa}\\
&\left[h+\mathcal{T}^h_b\right] u_t+gh(h+b)_x+huu_x+\mathcal{Q}^h u+\mathcal{Q}_b^h u=0\ , \label{eq:SGNb}
\end{align}
where 
\begin{align}
h(x,t)\doteq \eta(x,t)-b(x)\ , \label{eq:SGNc}
\end{align}
is the total depth of the water between the bottom $b(x)$ and the free surface elevation $\eta(x,t)$, $u(x,t)$ is the depth averaged horizontal velocity of the fluid, and $g$ the acceleration due to gravity. The operators $\mathcal{T}_b^h$, $\mathcal{Q}^h$ and $\mathcal{Q}_b^h$ depend on $h$, and are defined as follows:
\begin{equation}\label{eq:SGNd}
\mathcal{T}^h_b w= h\left[h_x b_x+\frac{1}{2} h b_{xx}+b_x^2\right]w-\frac{1}{3}\left[h^3 w_{x}\right]_x \ ,
\end{equation}
\begin{equation}\label{eq:SGNe} 
\mathcal{Q}^h w = -\frac{1}{3}\left[h^3(ww_{xx}-w_x^2)\right]_x\ ,
\end{equation}
\begin{equation}
 \mathcal{Q}_b^h w = \frac{1}{2}\left[h^2(w^2b_{xx}+ww_xb_x)\right]_x-\frac{1}{2}h^2\left[ww_{xx}-w_x^2\right]b_x + hw^2b_x b_{xx}+hww_xb_x^2\ . \label{eq:SGNf}
\end{equation}
\label{eq:Qboperator}
\end{subequations}
In dimensional and unscaled form, the independent variable, $x\in\mathbb{R}$ is a spatial variable and $t\geq 0$ represents the time.

The \acs{SGN} equations as derived by Seabra-Santos {\em et.al.} in \cite{SRT1987}  have also
been derived in a three-dimensional form in \cite{LB2009} and in a different formulation by Green and Naghdi \cite{GN1976}. In the case of a flat bottom (i.e. $b_x=0$) (\ref{eq:SGN})  is simplified to the so-called Serre system of equations derived first by Serre \cite{S1953I,S1953II} and re-derived later by Su and Gardner, \cite{Su1969}. For these reasons the equations (\ref{eq:SGN}) are also known as the Serre, or Green-Nagdhi,  or Su and Gardner equations. We will henceforth refer to them here as the \acf{SGN} equations.

Under the additional assumption of small amplitude waves (i.e. the solutions are of small amplitude), the \acs{SGN} system reduces to the Peregrine system, \cite{Peregrine1967}:
\begin{subequations} \label{eq:Peregrine}
\begin{align}
& h_t+(hu)_x =0\ , \label{eq:Peregrinea}\\
& u_t+g(h+b)_x+uu_x-\frac{b}{2}[bu]_{xxt}-\frac{b^2}{6}u_{xxt}=0\ . \label{eq:Peregrineb}
\end{align}
\end{subequations}
Peregrine's system belongs to the weakly dispersive and weakly nonlinear \acs{BT} systems. There are also other \acs{BT} systems that are asymptotically equivalent to Peregrine's system, cf. \cite{Mitsotakis2007,Nwogu1993}. Differences between the \acs{SGN} equations and \acs{BT} models are explained in \cite{DM2010}. In the same work the inclusion of surface tension effects have been included and explained in detail.

Although equations (\ref{eq:Peregrine}) can be derived from the \acs{SGN} system, their solutions  have different properties. For example the solutions of the \acs{SGN} equations are invariant under the Galilean boost, while the respective solutions of (\ref{eq:Peregrine}) are not. It is noted that the \acs{SGN} equations have a Hamiltonian formulation, \cite{Li2002,Israwi2011}. Specifically, for a stationary bathymetry the \acs{SGN} system conserves the total energy functional \cite{Israwi2011}:
\begin{equation}\label{eq:energy}
I(t)=\int_{\mathbb{R}}g\eta^2+hu^2+T_b^hu\cdot u~ dx\ ,
\end{equation}
in the sense that $I(t)=I(0)$, for all $t>0$. The conservation of this Hamiltonian will be used to measure  accuracy and conservation properties of the proposed numerical methods.

Although both systems are known to admit solitary wave solutions propagating without change in their shape over a horizontal bottom $y=-b_0$, only the solitary waves of the \acs{SGN} system have known formulas in a closed form. Specifically, a solitary wave of the \acs{SGN} system with amplitude $A$ can be written in the form:
\begin{subequations} \label{eq:sw}
\begin{align}
& h_s(x,t)=b_0+A{\sech}^2[\lambda(x-c_s t)],\  u_s(x,t)=c_0\left(1-\frac{b_0}{h(x,t)}\right)\ ,  \label{eq:swa}\\
&  \lambda = \sqrt{\frac{3A}{4b_0^2(b_0+A)}}\mbox{ and } c_s=c_0\sqrt{1+\frac{A}{b_0}}\ , \label{eq:swb}
\end{align}
\end{subequations}
where $c_s$ is the phase speed of the solitary wave and $c_0=\sqrt{gb_0}$ is the linear wave speed.
 
In addition to the above properties the system (\ref{eq:SGN}) is known to have several favourable well-posedness properties. For example, it is locally well-posed in time for smooth bottom functions in $\mathbb{R}$. Specifically, if $C_b^{\infty}(I)$ is the space of bounded and continuously differentiable functions on  the interval $I$ and $H^{s}(I)$ denotes the usual Sobolev space of $s$-order weakly differentiable functions on $I$, and if  $b\in C_b^{\infty}(\mathbb{R})$, $\kappa>1/2$, $s\geq \kappa+1$ with the initial condition is $(\eta_0,u_0)\in H^{s}(\mathbb{R})\times H^{s+1}(\mathbb{R})$, then there is at least a maximal time $T_{\max}>0$, such that the \acs{SGN} equations admit a unique solution $(\eta,u)\in   H^{s}(\mathbb{R})\times H^{s+1}(\mathbb{R})$, cf. \cite{Israwi2011}.
We also mention that the solution will satisfy a non-vanishing depth condition when the initial data satisfy the same non-vanishing condition, i.e. there is an $\alpha>0$ such that for all $t\geq 0$ 
\begin{equation}\label{depthcond}
\eta(x,t)-b(x)\geq \alpha>0\ .
\end{equation}  
For more information on the model equations, including the derivation, theory and justification, refer to  \cite{Lannes2013, LB2009, B2004}. 

We consider here the \acs{SGN} equations on a finite interval $I=[a,b]$. In addition to periodic boundary conditions, we will study physically relevant reflective boundary conditions. The linearized equations (\ref{eq:SGN}) about the trivial solution coincide with the linearized Peregrine system. Initial-boundary value problems for Peregrine's system have been studied in \cite{FP2005}, where the existence of solutions was proven for boundary data for $u(a,t)$ and $u(b,t)$. In a similar manner, the boundary conditions $u(a,t)=u(b,t)=0$ for $t\geq 0$ are sufficient for the system  (\ref{eq:SGN}) to model the reflection of the waves on solid wall boundaries. These specific wall boundary conditions have been used widely to describe reflection of waves on the boundaries for various numerical models including \acs{BT} systems such as the Nwogu system, \cite{KD2013}.

The numerical discretization of the \acs{SGN} equations is a challenging problem due to their complicated form. Not only do the operators in front of the temporal derivatives depend on the unknown function $h(x,t)$, but they also contain high order derivatives in nonlinear terms. Because of the presence of those high order spatial partial derivatives the solution to the system should be smooth. Recently several schemes have been proposed such as Finite Difference / Finite Volume Schemes \cite{Mirie1982,SRT1987,ACSSA1993, BCLMT2010, CBB1, CBB2, CLM2010} and Discontinuous Galerkin methods \cite{LGLX2014, PDZKWD2014, PDZKWD2014}. Although these methods are very useful to study practical problems, such as the runup of nonlinear waves on slopes, they can be highly dissipative methods in the sense that they introduce numerical dissipation or dispersion, due to the approximation of the nonlinear terms by dissipative flux functions and the use of low order slope limiters. These methods appear to have good conservation properties though when high-order WENO methods are used or other non-classical numerical fluxes are used, \cite{LaMa2015}.   

Some other highly accurate numerical methods have been developed for the Serre equations in the case of a horizontal bottom, such as spectral methods \cite{Dutykh2013} and Galerkin / \acf{FEM} \cite{MID2014}. Although these methods appear to have satisfactory conservation properties, it is very difficult to extend them to the \acs{SGN} equations especially in the case of two horizontal dimensions. For example the standard Galerkin method of \cite{MID2014} requires tensor products of cubic splines in order to be consistent in two spatial dimensions in a similar manner to \cite{DMS1}. 

In this study, a new modified Galerkin  / \acf{FEM} is proposed for the numerical solution of the \acs{SGN} equations. This method allows the use of low-order finite elements such as piecewise quadratic ($P^2$) or even piecewise linear ($P^1$) Lagrange finite elements. Two of the main advantages of this method is that it is highly accurate, and  it has very good conservation properties. Other advantages are the sparsity of the resulting linear systems, the low complexity of the algorithm  due to the use of low-order finite element spaces and finally its potential to be extended to the two-dimensional model equations of
\cite{Lannes2013}. 

Similar techniques that reduce the requirement of high-order finite elements (for example the use of cubic splines) have been used previously for weakly-nonlinear Boussinesq systems in \cite{WB1999,WB2002,DMS3} where the second derivative in the linear dispersive terms has been replaced by either the discrete Laplacian operator or the solution of an intermediate problem. In the case of the Bona-Smith type of Boussinesq systems with wall boundary conditions, the modified Galerkin method converges at an optimal rate showing great performance contrary to the suboptimal convergence rates achieved with the standard Galerkin / \acs{FEM} method, \cite{DMS3}. Similarly to the behavior of the modified Galerkin method for the weakly-nonlinear Boussinesq equations, the proposed \acs{FEM} scheme for the \acs{SGN} equations can achieve optimal convergence rates depending on the choice of the trial function spaces, contrary to the suboptimal rates obtained for the standard Galerkin method for the \acs{SGN} system and also for Peregrine's system as shown in \cite{AD2012}. Specifically, in order to achieve optimal convergence properties one may use spaces of piecewise linear elements for the free surface elevation and piecewise quadratic elements for the horizontal velocity.

The validated numerical method is applied to study the \acs{SGN} equations with reflective boundary conditions in a systematic way through a series of numerical experiments. In particular, we focus on the following issues:
\begin{itemize}
\item accuracy of the modified Galerkin method and invariant conservation;
\item reflection of solitary waves at a vertical wall; and
\item shoaling of solitary waves on plain or composite beaches.
\end{itemize}
The convergence properties of the new numerical method are also tested in the case of periodic boundary conditions. For more information about the behavior and the properties of the Galerkin / \acs{FEM} method with periodic boundary conditions we refer to \cite{MID2014}.

This paper is organized as follows. Section \ref{numethods} presents the fully discrete schemes for the \acs{SGN} equations. In Section \ref{accuandconv} we study the convergence, accuracy and numerical stability of the modified Galerkin method. Finally, Section \ref{sec:numexp} presents computational studies of shoaling and reflected waves validating both the choice of boundary conditions and the numerical scheme. We close the paper with conclusions in Section \ref{sec:summary}.

\section{The numerical methods}\label{numethods}

We consider the \acf{IBVP} comprising System (\ref{eq:SGN}), subject to reflective boundary conditions:
\begin{subequations} \label{eq:SGN2}
\begin{align}
& h_t+(hu)_x =0\ , \label{eq:SGN2a}\\
&\left[h+\mathcal{T}^h_b\right] u_t+gh(h+b)_x+huu_x+\mathcal{Q}^h u+\mathcal{Q}_b^h u=0\ , \label{eq:SGN2b}\\
& u(a,t)=u(b,t)=0, \  h(x,0)=h_0(x), \  u(x,0)=u_0(x)\ , \label{eq:SGN2c}
\end{align}
where again
\begin{equation}\label{eq:SGN2d}
\mathcal{T}^h_b w= h\left[h_x b_x+\frac{1}{2} h b_{xx}+b_x^2\right]w-\frac{1}{3}\left[h^3 w_{x}\right]_x \ ,
\end{equation}
\begin{equation}\label{eq:SGN2e} 
\mathcal{Q}^h w = -\frac{1}{3}\left[h^3(ww_{xx}-w_x^2)\right]_x\ ,
\end{equation}
\begin{equation}
 \mathcal{Q}_b^h w = \frac{1}{2}\left[h^2(w^2b_{xx}+ww_xb_x)\right]_x-\frac{1}{2}h^2\left[ww_{xx}-w_x^2\right]b_x + hw^2b_x b_{xx}+hww_xb_x^2\ , \label{eq:SGN2f}
\end{equation}
\label{eq:Qboperator}
\end{subequations}
$x\in(a,b)\subset \mathbb{R}$ and $t\in[0,T]$. We assume that (\ref{eq:SGN2}) possesses a unique solution, such that $h$ and $u$ are sufficiently smooth and, for any $t\in[0,T]$, in suitable Sobolev spaces
$$h(x,\cdot)\in H^s, \quad u(x,\cdot)\in H_0^{s+1}\ , $$
where $s\geq 1$. Here and below, $\|\cdot\|_s$ denotes the standard norm in $H^s$ while $H^s_0$  will denote the subspace of $H^s$ whose elements vanish at $x=a$ and $x=b$. We also use the inner product in $L^2\equiv H^0$, denoted by $(\cdot ,\cdot)$, which is
$$(u,v)=\int_a^b uv~dx\ ,$$
for all $u,v\in L^2$.

A spatial grid of the interval $[a,b]$ is a collection of points $x_i=a+i~ \Delta x$, for $i=0,1,\cdots,N$, where $\Delta x$ is the grid size, and $N\in \mathbb{N}$, such that $\Delta x=(b-a)/N$. Let $(\tilde{h},\tilde{u})\in S_h\times S_u$ be the corresponding spatially discretized solutions of the Galerkin / \acf{FEM} for suitable finite-dimensional spaces $S_h$ and $S_u$. First we present the standard Galerkin / \acs{FEM} semidiscretization. 

\subsection{The standard Galerkin method}\label{standard}
For the standard Galerkin method, we consider the space of smooth splines
$$S^r=\left\{ \left. \varphi\in C^{r-1}[a,b] \right| ~\varphi|_{[x_i,x_{i+1}]}\in \mathbb{P}^{r},~ 0\leq i\leq N-1\right\}\ ,$$ where $\mathbb{P}^r$ is the space of polynomials of degree $r$. The standard Galerkin method requires $r\geq 3$. Here, we take $r=3$. The trial function space for the first equation and for the solution $\tilde{h}$ is chosen as $S_h=S^r$, while for the second equation and for the approximation of the depth averaged velocity of the fluid $\tilde{u}$ is $S_u=S^r\cap\left\{\left. \varphi\in C[a,b]\right|~ \varphi(a)=\varphi(b)=0\right\}$. To state the associated semi-discrete problem, let $\phi\in S_h$ and $\psi\in S_u$ be arbitrary test functions. After taking inner products, and using integration by parts, the semi-discrete problem takes the form:
\begin{subequations} \label{eq:semi}
\begin{align}
& (\tilde{h}_t,\phi)+\left((\tilde{h}\tilde{u})_x,\phi \right) =0\ , \label{eq:semia}\\
& \mathcal{B}(\tilde{u}_t,\psi;\tilde{h})+\left(\tilde{h}\left[g(\tilde{h}+\tilde{b})_x+\tilde{u}\tilde{u}_x\right],\psi \right) 
+\mathcal{Q}(\tilde{u},\psi;\tilde{h})+\mathcal{Q}_b(\tilde{u},\psi;\tilde{h})=0\ , \label{eq:semib}
\end{align}
where $\mathcal{B}$, $\mathcal{Q}$ and $\mathcal{Q}_b$ are defined for any $\omega,\psi\in S_u$ as
\begin{align}
 \mathcal{B}(\omega,\psi;\tilde{h})\doteq & \left(\tilde{h}\left[1+\tilde{h}_x\tilde{b}_x+\frac{1}{2}\tilde{h}\tilde{b}_{xx}+\tilde{b}_x^2 \right]\omega,\psi \right)+\frac{1}{3}\left(\tilde{h}^3w_x,\psi_x \right)\ , \label{eq:semic}\\
 \mathcal{Q}(\omega,\psi;\tilde{h})\doteq & \frac{1}{3}\left(\tilde{h}^3\left[\omega\omega_{xx}-\omega_x^2 \right],\psi_x \right)\ , \label{eq:semid}\\
 \mathcal{Q}_b(\omega,\psi;\tilde{h})\doteq & -\frac{1}{2}\left(\tilde{h}^2 (\omega^2\tilde{b}_{xx}+\omega\omega_x\tilde{b}_x),\psi_x \right)- \nonumber \\& \frac{1}{2}\left(\tilde{h}\tilde{b}_x\left\{\tilde{h}\left[\omega\omega_{xx}-\omega_x^2\right]-\omega^2\tilde{b}_{xx}-\omega\omega_x\tilde{b}_x\right\},\psi\right)\ . \label{eq:semie}
\end{align}
\end{subequations}

This is a system of \acf{ODEs}. Given the initial conditions $\tilde{h}(x,0)\doteq \tilde{h}_0=\mathcal{P}_h(h_0)$ and $\tilde{u}(x,0)\doteq \tilde{u}_0=\mathcal{P}_u(u_0)$ where $\mathcal{P}_h$ and $\mathcal{P}_u$ are appropriate projections on $S_h$ and $S_u$ respectively, we assume that the system (\ref{eq:semi}) has a unique solution. Appropriate projections of the initial conditions could be the standard $L^2$-projections on $S_h$ and $S_u$, defined as $H_0\in S_h$ and $U_0\in S_u$ such that
$$\int_a^b H_0\phi~ dx = \int_a^bh_0\phi~dx, \mbox{ for all } \phi\in S_h\ ,$$
and 
$$\int_a^b U_0\psi~ dx = \int_a^bu_0\psi~dx, \mbox{ for all } \psi\in S_u\ .$$

The presence of the term $\omega_{xx}$  in $\mathcal{Q}$ and $\mathcal{Q}_b$ requires the use of at least $C^2$ smooth splines. In particular, one may use cubic splines, which correspond to $S_h \equiv S^3$, i.e., $r=3$. A basis for the space $S_h$  
for a uniform grid of mesh-length $\Delta x$ can be formed by the functions $\phi_j(x)=B(x-x_j/\Delta x)\Big|_{[a,b]}$,
 $j=-1,0,\cdots, N+1$ where $x_{-1}=a-\Delta x$, $x_{N+1}=b+\Delta x$ and 
\begin{equation*}
B(x)=\left\{ \begin{array}{lr}
\frac{1}{4}(x+2)^3, & -2\leq x\leq -1\ ,\\
\frac{1}{4}[1+3(x+1)+3(x+1)^2-3(x+1)^3], & -1\leq x\leq 0\ ,\\
\frac{1}{4}[1+3(1-x)+3(1-x)^2-3(1-x)^3], & 0\leq x\leq 1\ ,\\
\frac{1}{4}(2-x)^3, & 1\leq x\leq 2\ , \\
0, & x\in R-[-2,2]\ .
\end{array} \right.
\end{equation*}
The basis of $S_u=S^3_0$ can be described by the functions $\psi_j(x)=\phi_j(x)$ for $2\leq j\leq N-2$, plus four functions $\psi_0,\psi_1,\psi_{N-1}, \psi_{N}$,
 taken as linear combinations of the $\phi_{-1}, \phi_0,\phi_1$ and $\phi_{N-1}, \phi_N, \phi_{N+1}$, which are such that $\psi_0(a)=\psi_1(a)=\psi_{N-1}(b)=\psi_N(b)=0$.
 For example, we take $\psi_0=\phi_0-4\phi_{-1}$ and $\psi_1=\phi_1-\phi_{-1}$, \cite{Schultz1973}. The convergence properties of the standard Galerkin method with cubic splines are very similar to those of the modified Galerkin method with $S^3$ elements and thus are not presented here. Some details and numerical experiments with the standard Galerkin method can be found in \cite{KZM2015}.
 
\subsection{The modified Galerkin method}\label{modified}

In real-world applications, the use of low-order finite element methods can lead to faster computations. We derive a numerical method that does not require high-order finite element spaces. This can be done by using similar techniques to those proposed in \cite{WB1999,WB2002,DMS3}, but applied to the nonlinear term $uu_{xx}$. The resulting method is a modified Galerkin method that allows the use of Lagrange elements as low as $P^1$ in order, consisting of piecewise linear functions. We proceed with the derivation of the modified Galerkin method but first we introduce the notation for the Lagrange finite element spaces that are subspaces of $H^1$. It is noted that we cannot use Lagrange finite elements with the standard Galerkin method due to the second order derivative.

The space $P^r$ of Lagrange finite elements is defined on the grid $x_i=x_0+i\Delta x$, $i=0,1,\cdots, N$ as:
\begin{equation}\label{eq:Lagrange}
P^r=\left\{\chi\in C^0[a,b] \left. \right| \chi|_{[x_i,x_{i+1}]}\in \mathbb{P}^r,\ 0\le i\le N-1 \right\}\ ,
\end{equation}
We will restrict the analysis below to $P^1$, $P^2$ and $P^3$ Lagrange finite element spaces. The $P^1$ Lagrange elements can be defined using the basis functions
$$\chi_i(x)=\left\{ \begin{array}{l}
\frac{x-x_{i-1}}{\Delta x}, \ \mbox{ if } x\in [x_{i-1},x_i]\ ,\\
\frac{x_{i+1}-x}{\Delta x}, \ \mbox{ if } x\in [x_i, x_{i+1}]\ ,\\
0 \ , \mbox{ otherwise } \ .
\end{array} \right.$$ 
The $P^2$ Lagrange finite element space can be defined on the same grid $x_i=x_0+i\Delta x$  for $i=0,1,\cdots, N$ and at the midpoints $x_{i+1/2}=x_i+\Delta x/2$ for $i=0,1,\cdots, N-1$ using the basis functions 
$$ \chi_i(x)=\varphi\left(\frac{x-x_i}{\Delta x} \right), \ i=0,1,\cdots, N\ ,$$ and
$$\chi_{i+1/2}(x)=\psi\left(\frac{x-x_{i+1/2}}{\Delta x} \right),\ i=0,1,\cdots, N-1 \ ,$$
with
$$\varphi(x)=\left\{
\begin{array}{ll}
(1+x)(1+2x), & -1\leq x\leq 0\ ,\\
(1-x)(1-2x), & 0\leq x\leq 1\ ,\\
0& \mbox{otherwise}\ ,
\end{array} \right. \, $$
and 
$$\psi(x)=\left\{
\begin{array}{ll}
1-4x^2, & |x|\leq 1/2\ ,\\
0 , & \mbox{otherwise}\ .
\end{array} \right. $$ Then a function $w \in P^2$ is written as
$$w(x)=\sum_{i=0}^{N} w(x_i)\chi_i(x)+\sum_{i=0}^{N-1} w(x_{i+1/2})\chi_{i+1/2}(x)\ .$$
For the construction of the Lagrange basis function of the general space $P^r$ we refer to  \cite{EGbook}. We also consider the Lagrange finite element spaces $P^r_0=\{\chi \in P^r| \chi(a)=\chi(b)=0\}$. These spaces are subspaces of $H^1_0$ and will be used to approximate the depth averaged horizontal velocity of the water. For example, we will use the spaces $S_h=P^r$ and $S_u=P^q_0$, for some integers $r,q$ or we will use the spaces of smooth splines described in the previous section.
For more information related to Lagrange finite element spaces and its approximation properties, we refer to \cite{EGbook}. It is worth mentioning that, for the discretization of the bottom boundary, we use the $L^2$ projection of the bathymetry and its derivatives.

Given $v \in H^1_0$, we define the {\em non-linear discrete Laplacian} operator $\tilde{\partial}^2 :H_0^1\rightarrow S_u$ such that
\begin{equation}\label{eq:dlaplacian}
\left(\tilde{\partial}^2 v,\psi\right)= -\left(v_x^2,\psi\right)-\left(vv_x,\psi_x\right), \mbox{ for all }\psi\in S_u\ .
\end{equation}
In fact, the function $\tilde{\partial}^2v \in S_u$ approximates the function $vv_{xx}$ as if $v$ was a smooth $C^2$ function and substitution in (\ref{eq:semid}) and (\ref{eq:semie}) leads to the modified Galerkin semi-discretization:
\begin{subequations} \label{eq:semi2}
\begin{align}
& (\tilde{h}_t,\phi)+\left((\tilde{h}\tilde{u})_x,\phi \right) =0, \quad \phi\in S_h \ , \label{eq:semia2}\\
& \mathcal{B}(\tilde{u}_t,\psi;\tilde{h})+\left(\tilde{h}\left[g(\tilde{h}+\tilde{b})_x+\tilde{u}\tilde{u}_x\right],\psi \right) 
+ \tilde{\mathcal{Q}}(\tilde{u},\psi;\tilde{h})+ \tilde{\mathcal{Q}}_b(\tilde{u},\psi;\tilde{h})=0,\quad \psi\in S_u \ , \label{eq:semib2}
\end{align}
where $\mathcal{B}$, $\tilde{\mathcal{Q}}$ and $\tilde{\mathcal{Q}}_b$ are defined for any $\omega,\psi\in S_u$ as
\begin{equation} \label{eq:semic2}
 \mathcal{B}(\omega,\psi;\tilde{h})\doteq \left(\tilde{h}\left[1+\tilde{h}_x\tilde{b}_x+\frac{1}{2}\tilde{h}\tilde{b}_{xx}+\tilde{b}_x^2 \right]\omega,\psi \right) +\frac{1}{3}\left(\tilde{h}^3w_x,\psi_x \right) \ ,
 \end{equation}
 \begin{equation} \label{eq:semid2}
 \tilde{\mathcal{Q}}(\omega,\psi;\tilde{h})\doteq  \frac{1}{3}\left(\tilde{h}^3\left[\tilde{\partial}^2\omega-\omega_x^2 \right],\psi_x \right)\ , 
 \end{equation}
 \begin{align}
  \tilde{\mathcal{Q}}_b(\omega,\psi;\tilde{h})\doteq & -\frac{1}{2}\left(\tilde{h}^2 (\omega^2\tilde{b}_{xx}+\omega\omega_x\tilde{b}_x),\psi_x \right)- \nonumber \\ & \frac{1}{2}\left(\tilde{h}\tilde{b}_x\left\{\tilde{h}\left[\tilde{\partial}^2\omega-\omega_x^2\right]-\omega^2\tilde{b}_{xx}-\omega\omega_x\tilde{b}_x\right\},\psi\right)\ . \label{eq:semie2}
\end{align}
\end{subequations}
Although we assume that the bottom function has the appropriate smoothness, the absence of second order spatial derivatives of the depth integrated horizontal velocity in the semidiscrete scheme allows the use of Lagrange elements, as well as high-order elements such as cubic or quintic splines. 
In what follows, we test the efficiency of the modified Galerkin method using the spaces $P^1$, $P^2$ and $P^3$ of Lagrange elements and the space $S^3$ of cubic splines with periodic and reflective boundary conditions. 

\begin{remark}[Mass lumping]
In order to compute the non-linear discrete Laplacian $\tilde{\partial} w$ one needs to solve the linear system obtained by the discretization of (\ref{eq:dlaplacian}). In the case of wall boundary conditions if, for example, $\psi_i$ denotes basis functions of $S_u$ the system can be written as $\mathcal{M}w=f$ where the mass matrix is a banded matrix with entries $\mathcal{M}_{ij}=(\psi_i,\psi_j)$ and $f_i=-(v_x^2,\psi_j)-(vv_x,\psi_j')$. To improve the speed of the numerical method, one may apply the method of mass lumping in the formulation of the matrix $\mathcal{M}$. This can be done, for example in the case of quadratic Lagrange elements, by approximating the integrals with Simpson's rule. This leads to a diagonal matrix that can be inverted trivially. All the numerical experiments with $P^1$ and $P^2$ elements have been performed with mass lumping, in addition to the standard matrix formulation with comparable results. 
\end{remark}

\begin{remark}
The choice of the discrete Laplacian is not unique.  For example, when periodic boundary conditions are used, then one may consider the linear discrete Laplacian, which replaces the term $u_{xx}$ as proposed in \cite{WB1999,WB2002,DMS3}.
\end{remark}

\begin{remark}
The presence of the term $b_{xx}$ in the semidiscrete scheme implies the typical requirement of a smooth bottom. When a piecewise linear bottom topography is given, then the use of an appropriate projection onto the finite element space (such as the elliptic projection) or local smoothing of the bottom is required. In this paper the bathymetry is usually a piecewise linear function and thus we use the smoothing method described in \cite{AA2006}. 
\end{remark}

\begin{remark}
The modified Galerkin method can be expressed in an equivalent way by introducing an additional independent variable $w$ satisfying the equation $w=uu_{xx}$, \cite{WB1999}. This increase of the degrees of freedom in the model equations by one allows the approximation of the unknown function $u$ by low-order finite element spaces such as $P^1$ or $P^2$ finite element spaces.
\end{remark}

\subsection{Temporal discretization}

In \cite{MID2014}, it was shown by numerical means that the standard Galerkin method for the discretization of the \acs{SGN} equations with flat bottom leads to a system of \acs{ODEs} that is not stiff. Also, the classical, explicit, four-stage, fourth order Runge-Kutta (RK) method, described by the following {\em Butcher tableau}:
\begin{equation}
\label{tab:RK}
\begin{tabular}{c | c}
  $A$ & $b$ \\ \hline
  $\tau$ &  
\end{tabular} \ = \
\begin{tabular}{c c c c| c}
  0 & 0 & 0 & 0 & 1/6\\
  1/2 & 0 & 0 & 0 & 1/3 \\ 
  0 & 1/2 & 0 & 0 & 1/3 \\
  0 & 0 & 1/2 & 0 &1/6\\\hline
  0 & 1/2 & 1/2 & 1 &
\end{tabular}~,
\end{equation}
is able to integrate the respective semi-discrete system numerically in time without imposing restrictive stability conditions on the ratio $\Delta t/\Delta x$, but only mild restrictions on the mesh length such as $\Delta t/\Delta x\leq 2$ for smooth solutions.

Concerning the time integration, the modified Galerkin method has very similar behavior to the standard Galerkin method. Upon choosing appropriate basis functions for the spaces $S_h$ and $S_u$, the semidiscrete system (\ref{eq:semi2}) represents a system of ODEs. We use a uniform time-step $\Delta t$ such that $\Delta t=T/K$ for $K\in \mathbb{N}$. The temporal grid is then $t^n=n\Delta t$, where $n=0,1,\cdots, K$. Given the ODE $y'=\Phi(t,y)$, one step of this four-stage RK scheme (with $y^n$ approximating $y(t^n)$) is:
\begin{algorithmic}
  \For{$i = 1 \to 4$} 
   \State $\tilde{y}^i ~~~=\ y^n+\Delta t\sum_{j=1}^{i-1}a_{ij}\,y^{n,j}$ \vspace{1mm}
   \State $y^{n,i}  \ =\ \Phi(t^{n,i},\tilde{y}^i),  \quad \mbox{evaluated at~}  t^{n,i} \equiv t^n+\tau_i \Delta t $
  \EndFor 
  \State $y^{n+1}=y^n+\Delta t\sum_{j=1}^4 b_j\,y^{n,j}$~,  \vspace{2mm}
\end{algorithmic}
where $a_{ij}$, $\tau_i$, $b_i$ are given in table~(\ref{tab:RK}). Applying this scheme to~(\ref{eq:semi2}) 
and denoting by $H^n$ and $U^n$ 
the fully discrete approximations in $S_h$ and $S_u$ of 
$h(\cdot,t^n)$, $u(\cdot,t^n)$, respectively, 
leads to Algorithm~\ref{alg:Serre}.
\begin{algorithm}
\caption{Time-marching \acs{FEM} scheme for the IBVP of the system~\eqref{eq:semi2}}
\label{alg:Serre}
\begin{algorithmic}[h]
  \State $H^0 \ =\  \P \{h_0\}$ \vspace{2mm}
  \State $U^0 ~\, =\  \P \{u_0\}$ \vspace{2mm}
  \For{$n = 0 \to N-1$}  \vspace{2mm}
  \For{$i = 1 \to 4$}  \vspace{2mm}
    \State $\tilde{H}^i \ =\ H^n \ +\ \Delta t\sum_{j=1}^{i-1}a_{ij}\,H^{n,j}$\vspace{2mm}
    \State $\tilde{U}^i ~\, =\ U^n \ + \ \Delta t\sum_{j=1}^{i-1}a_{ij}\,U^{n,j}$\vspace{2mm}
    \State $\left(\tilde{\partial}^2U^{i},\phi \right)=-((U^{i}_x)^2,\phi)-(U^{i}U^{i}_x,\phi_x)$\vspace{2mm}
    \State $(H^{n,i},\psi) \hspace{5.5mm} \ =\ -((\tilde{H}^{i}\tilde{U}^{i})_x,\psi), 
    \quad\mbox{evaluated at~} t^{n,i} \equiv t^n+\tau_i \Delta t $\vspace{2mm}
    \State $\mathcal{B}(U^{n,i},\phi;\tilde{H}^{i}) \ =\ -\left(\tilde{H}^{i}\left[g(\tilde{H}^{i}-\tilde{b})_x+\tilde{U}^{i}\tilde{U}^{i}_x\right],\phi\right) 
    \ -\ \tilde{\mathcal{P}}(\tilde{U}^{i},\phi;\tilde{H}^{i},\tilde{b})=0$\vspace{2mm}
    \State where $\tilde{\mathcal{P}}(\tilde{U}^{i},\phi;\tilde{H}^{i},\tilde{b})=\tilde{\mathcal{Q}}(\tilde{U}^{i},\phi;\tilde{H}^{i})+\tilde{\mathcal{Q}}_{\tilde{b}}(\tilde{U}^{i},\phi;\tilde{H}^{i})$ \vspace{2mm}
  \EndFor \vspace{2mm}
  \State $H^{n+1} \ =\ H^n \ + \ \Delta t\sum_{j=1}^4 b_j\,H^{n,j}$ \vspace{2mm}
  \State $U^{n+1} ~\, =\ U^n \ +  \ \Delta t\sum_{j=1}^4 b_j\,U^{n,j}~$\vspace{2mm}
  \EndFor 
\end{algorithmic}
\end{algorithm}
Given bases $\{\phi_i\}$ of $S_h$ and $\{\psi_i\}$ of $S_u$, the implementation of Algorithm~\ref{alg:Serre} 
requires solving at each time step the following linear systems:
\begin{enumerate}
  \item[(a)] Four linear systems with the time-independent matrix $(\phi_i,\phi_j)$;
  \item[(b)] Four linear systems with the time-dependent matrix $\mathcal{B}(\psi_i,\psi_j;h)$;
  \item[(c)] Four linear systems with the time-independent matrix $(\psi_i,\psi_j)$~.
\end{enumerate}
The four linear systems in (a) arise from the discretization of the equation (\ref{eq:semia2}), while the four linear systems in (c) arise from the computation of the discrete Laplacian (\ref{eq:dlaplacian}). All of these matrices are banded and symmetric, and one can use either  direct methods for banded systems or classical iterative methods for sparse systems. To approximate the inner products, we use Gauss-Legendre quadrature 
with 3 nodes per $\Delta x$ for $P^1$ elements, 5 nodes for $P^2$ and $P^3$, and 8 nodes for $S^3$ elements. It is noted that most of the experiments with $P^1$, $P^2$ and mixed elements have also been tested using adaptive time-stepping methods such as the Runge-Kutta-Fehlberg method, \cite{Hairer2009}, to ensure that the errors introduced by the temporal integration are negligible.

\section{Accuracy and convergence}\label{accuandconv}

As it was pointed out in the introduction, the \acf{cB} system is similar in structure to the \acs{SGN} system, given that both systems admit the same number of boundary conditions, and the equation for the free surface is the same. We thus expect similar behavior for the convergence of the numerical method. In \cite{ADM3}, it was shown that the convergence of the Galerkin / \acs{FEM} method for the \acs{cB} system with periodic boundary conditions is optimal, while in \cite{AD2013,AD2012} it was shown that the convergence of the Galerkin / \acs{FEM} method with cubic splines for the \acs{cB} system is suboptimal, and that in order to achieve an optimal rate of convergence, a non-standard method should be used. Specifically, it was shown that for wall boundary conditions in the case of $P^1$ elements, the following error estimates hold:
\begin{equation}
\max_{0\leq t\leq T} \|h-\tilde{h}\|\leq C~\Delta x^{3/2}, \quad \max_{0\leq t\leq T}\|u-\tilde{u}\|\leq C~\Delta x^2\ ,
\end{equation}
 which are suboptimal for $h$ and optimal for $u$. In the case of cubic splines, the following error estimates hold: 
 \begin{equation}
 \max_{0\leq\leq T} \|h-\tilde{h}\|\leq C~\Delta x^{3.5}\sqrt{\ln \frac{1}{\Delta x}}, \quad \max_{0\leq t\leq T}\|u-\tilde{u}\|\leq C~\Delta x^{4}\sqrt{\ln \frac{1}{\Delta x}} \ ,
 \end{equation} which is suboptimal in both $h$ and $u$ but the factor $\ln \Delta x$ is not dominant and generally it has an effect on the accuracy of the method that is too small to observe. In \cite{KZM2015}, it is shown, by numerical means that the standard Galerkin method with cubic splines for the \acs{SGN} system with wall boundary conditions appears to have similar convergence properties.
 
The standard Galerkin method for the initial-periodic boundary value problem for the \acs{SGN} system has been studied in \cite{MID2014}, where it was shown that the convergence is optimal. It is worth mentioning that the approximation (\ref{eq:dlaplacian}) can be used also with periodic boundary conditions, and the resulting modified Galerkin method has optimal convergence rates including all the advantages of the modification. Specifically, the errors in the $L^2$ norm for the periodic case are of $O(\Delta x^2)$ for the $P^1$ elements for both $h$ and $u$, $O(\Delta x^3)$ for the $P^2$ elements and for the $S^3$ elements is of $O(\Delta x^4)$. Because the results for the periodic problem are very similar to those of \cite{MID2014} and to the modified Galerkin method for the \acs{SGN} system with periodic boundary conditions, we don't present them here. 

We continue with the wall boundary conditions. In order to study the convergence of the numerical method, we use the \acs{SGN} equations in nondimensional but unscaled form. We start with the \acs{IBVP} (\ref{eq:SGN2}). Because we don't know any analytical solutions for this problem, we  consider the non-homogenous problem, which has the exact solution $h(x,t)= 1+e^{2t}(\cos(\pi x)+x+2)$ and $u(x,t)=e^{-tx}x\sin(\pi x)$ satisfying the equations (\ref{eq:SGN2a}) and (\ref{eq:SGN2b}) with the appropriate right-hand sides. The solution is computed using the modified Galerkin method in the interval $[0,1]$ and for $t\in (0,T]$ with $T=1$. The  finite element spaces of our choice are $P^1$ and $P^2$ Lagrange finite element spaces and the space $S^3$ of cubic splines. It is noted that we also tried other ``exact'' solutions with different right-hand sides in (\ref{eq:SGN2a}) and (\ref{eq:SGN2b}) to verify the computed rates of convergence. The results have been all very similar, and thus we only present the results of one case.

In order to study the accuracy and the convergence of the modified Galerkin method, several error indicators have been computed. The computed errors are normalized and defined as
\begin{equation}\label{eq:errsnorm}
E_s[F]\doteq \frac{\|F(x,T;\Delta x)-F_{\mbox{exact}}(x,T)\|_s}{\|F_{\mbox{exact}}(x,T)\|_s}\ ,
\end{equation}
where $F=F(\cdot ;\Delta x)$ is the computed solution, i.e., either $H\approx h(x,T)$ or $U\approx u(x,T)$, $F_{\mbox{exact}}$ is the corresponding exact solution and $s=0,1,2,\infty$ correspond to the $L^2$, $H^1$, $H^2$ and $L^{\infty}$ norms, respectively.  The analogous rates of convergence are defined as
\begin{equation}
\mbox{rate for }E_s[F]\doteq \frac{\ln (E_s[F(\cdot;\Delta x_{k-1})]/E_s[F(\cdot;\Delta x_k)])}{\ln(\Delta x_{k-1}/\Delta x_k)}\ ,
\end{equation} 
where $\Delta x_k$ is the grid size listed in row $k$ in each table. To ensure that the errors incurred by the temporal integration do not affect the rates of convergence we use $\Delta t\ll \Delta x$ while we take $\Delta x=1/N$. 

Tables \ref{tab:perconvlin1}, \ref{tab:perconvlin2} and \ref{tab:perconvlin3} present the errors and the corresponding rates of convergence of the modified Galerkin method with $P^1$ finite elements. It is shown that the rate of convergence is suboptimal for the total depth $h$ and optimal for velocity $u$. Specifically, Table \ref{tab:perconvlin1} suggests that $\|h-\tilde h\|=O(\Delta x ^{3/2})$ and $\|u-\tilde{u}\|=O(\Delta x^2)$. Table \ref{tab:perconvlin2} suggests that $\|h-\tilde{h}\|_1=O(\Delta x^{1/2})$ and $\|u-\tilde{u}\|_1=O(\Delta x)$. Finally, the $L^\infty$ estimates shown in \ref{tab:perconvlin3} are worse than the $L^2$ estimates, suggesting $\|h-\tilde{h}\|_\infty=O(\Delta x)$ and $\|u-\tilde{u}\|_\infty=O(\Delta x^2)$. All these results are very similar to those obtained theoretically and numerically in the case of Peregrine's system \cite{AD2013}.

\begin{table}[ht]
\centering
\begin{tabular}{llclc}
\hline
$N$ & $E_0[H]$ & rate for $E_0[H]$ & $E_0[U]$ & rate for $E_0[U]$\\
\hline
$10$ & $1.4661\times 10^{-2}$ & -- & $2.9141\times 10^{-2}$ & -- \\
$20$ & $3.3761\times 10^{-3}$  & $2.1186$ & $2.3778\times 10^{-3}$ & $3.6154$ \\
$40$ & $1.1967\times 10^{-3}$  & $1.4963$ & $5.6477\times 10^{-4}$ & $2.0739$ \\ 
$80$ & $4.4536\times 10^{-4}$  & $1.4260$ & $1.3856\times 10^{-4}$ & $2.0271$ \\
$160$ & $1.6179\times 10^{-4}$ &  $1.4609$ & $3.4262\times 10^{-5}$ & $2.0159$ \\ 
$320$ & $5.7982\times 10^{-5}$ &  $1.4804$ & $8.5165\times 10^{-6}$ & $2.0083$ \\ 
$640$ & $2.0638\times 10^{-5}$ &  $1.4903$ & $2.1229\times 10^{-6}$ & $2.0042$ \\ 
\hline 
\end{tabular}
\caption{Spatial errors and rates of convergence for the exact solution with $P^1$ finite elements using the $L^2$ norm.} \label{tab:perconvlin1}
\end{table}
\begin{table}[ht]
\centering
\begin{tabular}{llclc}
\hline
$N$ & $E_1[H]$ & rate for $E_1[H]$ & $E_1[U]$ & rate for $E_1[U]$\\
\hline
$10$ & $2.8836\times 10^{-1}$ & -- & $1.6945\times 10^{-1}$ & -- \\ 
$20$ & $1.8855\times 10^{-1}$ & $0.6129$ & $6.6788\times 10^{-2}$ & $1.3433$  \\ 
$40$ & $1.3684\times 10^{-1}$ & $0.4625$ & $3.2009\times 10^{-2}$ & $1.0611$  \\ 
$80$ & $1.0135\times 10^{-1}$ & $0.4331$ & $1.5859\times 10^{-2}$ & $1.0131$  \\ 
$160$ & $7.3388\times 10^{-2}$ & $0.4658$ & $7.9074\times 10^{-3}$ & $1.0041$  \\ 
$320$ & $5.2503\times 10^{-2}$ & $0.4831$ & $3.9492\times 10^{-3}$ & $1.0016$  \\ 
$640$ & $3.7340\times 10^{-2}$ & $0.4917$ & $1.9736\times 10^{-3}$ & $1.0007$ \\ 
\hline 
\end{tabular}
\caption{Spatial errors and rates of convergence for the exact solution with $P^1$ finite elements using the $H^1$ norm.}\label{tab:perconvlin2}
\end{table}
\begin{table}[ht]
\centering
\begin{tabular}{llclc}
\hline
$N$ & $E_\infty[H]$ & rate for $E_\infty[H]$ & $E_\infty[U]$ & rate for $E_\infty[U]$\\
\hline
$10$ & $3.4775\times 10^{-2}$ & $1.4587$ & $3.9193\times 10^{-2}$ & $1.4068$ \\
$20$ & $1.1039\times 10^{-2}$ & $1.6554$ & $5.7364\times 10^{-3}$ & $2.7724$ \\   
$40$ & $5.8666\times 10^{-3}$ & $0.9121$ & $1.0352\times 10^{-3}$ & $2.4701$ \\   
$80$ & $3.4929\times 10^{-3}$ & $0.7481$ & $2.1191\times 10^{-4}$ & $2.2884$ \\   
$160$ & $1.8861\times 10^{-3}$ & $0.8890$ & $5.2158\times 10^{-5}$ & $2.0225$ \\   
$320$ & $9.7698\times 10^{-4}$ & $0.9490$ & $1.3480\times 10^{-5}$ & $1.9520$ \\  
$640$ & $4.9684\times 10^{-4}$ & $0.9755$ & $3.4245\times 10^{-6}$ & $1.9769$ \\ 
\hline 
\end{tabular}
\caption{Spatial errors and rates of convergence for the exact solution with $P^1$ finite elements using the $L^\infty$ norm.}\label{tab:perconvlin3}
\end{table}

In the case of quadratic $P^2$ elements, the results with the same initial conditions are similar, but the convergence rates for $h$ are different to the expected suboptimal rates, analogous to the $P^1$ case. The respective convergence rates for $u$ are optimal as expected. Tables \ref{tab:perconvquad1}, \ref{tab:perconvquad2} and \ref{tab:perconvquad3} present the errors and the convergence rates for the modified Galerkin method with $P^2$ elements. Table \ref{tab:perconvquad1} suggests that the $\|h-\tilde h\|=O(\Delta x^2)$ and $\|u-\tilde{u}\|=O(\Delta x^3)$. Table \ref{tab:perconvquad2} suggests that $\|h-\tilde{h}\|_1=O(\Delta x)$ and $\|u-\tilde{u}\|=O(\Delta x^2)$. Finally, Table \ref{tab:perconvquad3} suggests the same estimates in the $L^\infty$ norm with $L^2$ norm, i.e. $\|h-\tilde{h}\|_\infty=O(\Delta x^2)$ and $\|u-\tilde{u}\|_\infty=O(\Delta x^3)$.  The same behaviour was observed when we used $P^3$ elements, i.e. the convergence rate for $h$ was suboptimal $\|h-\tilde h\|_s=O(\Delta x^{3-s})$ while the convergence for $u$ was optimal $\|u-\tilde{u}\|_s=O(\Delta x^{4-s})$, for $s=0,1$.
\begin{table}[ht]
\centering
\begin{tabular}{llclc}
\hline
$N$ & $E_0[H]$ & rate for $E_0[H]$ & $E_0[U]$ & rate for $E_0[U]$\\
\hline
$10$ & $3.1489\times 10^{-3}$ & -- & $1.5559\times 10^{-3}$ & -- \\   
$20$ & $6.1429\times 10^{-4}$ & $2.3579$ & $1.1718\times 10^{-4}$ & $3.7310$ \\   
$40$ & $1.3897\times 10^{-4}$ & $2.1441$& $9.8909\times 10^{-6}$ & $3.5665$ \\   
$80$ & $3.3909\times 10^{-5}$ & $2.0351$ & $1.0744\times 10^{-6}$ & $3.2025$ \\  
$160$ & $8.4285\times 10^{-6}$ & $2.0083$ & $1.2858\times 10^{-7}$ & $3.0628$ \\  
$320$ & $2.1019\times 10^{-6}$ & $2.0035$ & $1.5874\times 10^{-8}$ & $3.0180$ \\ 
$640$ & $5.2473\times 10^{-7}$ & $2.0021$ & $1.9773\times 10^{-9}$ & $3.0051$ \\ 
\hline 
\end{tabular}
\caption{Spatial errors and rates of convergence for the exact solution with $P^2$ finite elements using the $L^2$ norm.}\label{tab:perconvquad1}
\end{table}
\begin{table}[ht]
\centering
\begin{tabular}{llclc}
\hline
$N$ & $E_1[H]$ & rate for $E_1[H]$ & $E_1[U]$ & rate for $E_1[U]$\\
\hline
$10$ & $1.6038\times 10^{-1}$ &  -- & $1.8439\times 10^{-2}$ & --   \\ 
$20$ & $7.6716\times 10^{-2}$ & $1.0640$ & $3.3234\times 10^{-3}$ & $2.4721$   \\ 
$40$ & $3.7406\times 10^{-2}$ &  $1.0362$ & $6.4326\times 10^{-4}$ & $2.3692$   \\ 
$80$ & $1.8545\times 10^{-2}$ &  $1.0122$ & $1.4083\times 10^{-4}$ & $2.1915$   \\ 
$160$ & $9.2498\times 10^{-3}$ &  $1.0036$ & $3.3536\times 10^{-5}$ & $2.0701$  \\ 
$320$ & $4.6179\times 10^{-3}$ &  $1.0022$ & $8.2655\times 10^{-6}$ & $2.0206$ \\ 
$640$ & $2.3064\times 10^{-3}$ &  $1.0016$ & $2.0585\times 10^{-6}$ & $2.0055$  \\ 
\hline 
\end{tabular}
\caption{Spatial errors and rates of convergence for the exact solution with $P^2$ finite elements using the $H^1$ norm.}\label{tab:perconvquad2}
\end{table}
\begin{table}[ht]
\centering
\begin{tabular}{llclc}
\hline
$N$ & $E_\infty[H]$ & rate for $E_\infty[H]$ & $E_\infty[U]$ & rate for $E_\infty[U]$\\
\hline
$10$ & $7.4168\times 10^{-3}$ & -- & $2.7712\times 10^{-3}$ &  -- \\ 
$20$ & $1.7953\times 10^{-3}$ & $2.0466$ & $2.4407\times 10^{-4}$ & $3.5051$\\
$40$ & $4.0253\times 10^{-4}$ & $2.1571$ & $2.3644\times 10^{-5}$ & $3.3678$ \\  
$80$ & $1.0353\times 10^{-4}$ & $1.9589$ & $1.8056\times 10^{-6}$ & $3.7109$ \\  
$160$ & $2.5922\times 10^{-5}$ & $1.9979$ & $1.6493\times 10^{-7}$ & $3.4525$ \\ 
$320$ & $6.5561\times 10^{-6}$ & $1.9833$ & $2.0640\times 10^{-8}$ & $2.9984$ \\ 
$640$ & $1.6390\times 10^{-6}$ & $2.0000$ & $2.5810\times 10^{-9}$ & $2.9994$ \\ 
\hline 
\end{tabular}
\caption{Spatial errors and rates of convergence for the exact solution with $P^2$ finite elements using the $L^\infty$ norm.}\label{tab:perconvquad3}
\end{table}

We also studied the errors and convergence rates of a mixed modified Galerkin method, where we considered $S_h=P^1$ and $S_u=P^2$ i.e. we used linear elements for the approximation of $h$ and quadratic elements for the approximation of $u$. This was the only case where optimal rates of convergence for both $h$ and $u$ were obtained. Tables \ref{tab:perconvmix1}--\ref{tab:perconvmix3} suggest that $\|h-\tilde{h}\|_s=O(\Delta x^{2-s})$ and $\|u-\tilde{u}\|_s=O(\Delta x^{3-s})$ for $s=0,1$ while in the maximum norm we found that $\|h-\tilde{h}\|_\infty =O(\Delta x^2)$ and $\|u-\tilde{u}\|_\infty=O(\Delta x^3)$. Using mixed elements of higher order does not give analogously optimal results. For example, when we took $S_h=P^2$ and $S_u=P^3$ the convergence was suboptimal for $h$ and optimal for $u$. More specifically, the errors were  $\|h-\tilde{h}\|_s=O(\Delta x^{2-s})$ and $\|u-\tilde{u}\|_s=O(\Delta x^{4-s})$ for $s=0,1$ while in the maximum norm we found that $\|h-\tilde{h}\|_\infty =O(\Delta x^2)$ and $\|u-\tilde{u}\|_\infty=O(\Delta x^4)$. Observe that the convergence for $h$ achieved with the mixed $P^1-P^2$ elements  was the same with the mixed $P^2-P^3$ elements.
\begin{table}[ht]
\centering
\begin{tabular}{llclc}
\hline
$N$ & $E_0[H]$ & rate for $E_0[H]$ & $E_0[U]$ & rate for $E_0[U]$\\
\hline
$10$ & $1.3308\times 10^{-3}$ & -- & $1.4064\times 10^{-3}$ & -- \\ 
$20$ & $2.6812\times 10^{-4}$ & $2.3113$ & $9.9809\times 10^{-5}$ & $3.8167$ \\
$40$ & $6.2091\times 10^{-5}$ & $2.1104$ & $9.1910\times 10^{-6}$ & $3.4409$ \\ 
$80$ & $1.5350\times 10^{-5}$ & $2.0161$ & $1.0541\times 10^{-6}$ & $3.1241$ \\
$160$ & $3.8219\times 10^{-6}$ & $2.0059$ & $1.2825\times 10^{-7}$ & $3.0390$ \\ 
$320$ & $9.5358\times 10^{-7}$ & $2.0029$ & $1.6047\times 10^{-8}$ & $2.9986$ \\ 
$640$ & $2.38157\times 10^{-7}$ & $2.0014$ & $1.9751\times 10^{-9}$ & $3.0223$ \\ 
\hline 
\end{tabular}
\caption{Spatial errors and rates of convergence for the exact solution with $S_h=P^1$ and $S_u=P^2$ using the $L^2$ norm.}\label{tab:perconvmix1}
\end{table}
\begin{table}[ht]
\centering
\begin{tabular}{llclc}
\hline
$N$ & $E_1[H]$ & rate for $E_1[H]$ & $E_1[U]$ & rate for $E_1[U]$\\
\hline
$10$ & $6.9247\times 10^{-2}$ & -- & $1.5110\times 10^{-2}$ & -- \\ 
$20$ & $3.3950\times 10^{-2}$ & $1.0283$ & $2.6966\times 10^{-3}$ & $2.4863$ \\ 
$40$ & $1.6779\times 10^{-2}$ & $1.0168$ & $5.7095\times 10^{-4}$ & $2.2397$ \\   
$80$ & $8.3767\times 10^{-3}$ & $1.0022$ & $1.3474\times 10^{-4}$ & $2.0831$  \\  
$160$ & $4.1860\times 10^{-3}$ & $1.0008$ & $3.3114\times 10^{-5}$ & $2.0247$ \\  
$320$ & $2.0924\times 10^{-3}$ & $1.0004$ & $8.3275\times 10^{-6}$ & $1.9915$  \\ 
$640$ & $1.0460\times 10^{-3}$ & $1.0002$ & $2.0570\times 10^{-6}$ & $2.0173$ \\  
\hline 
\end{tabular}
\caption{Spatial errors and rates of convergence for the exact solution with $S_h=P^1$ and $S_u=P^2$ using the $H^1$ norm.}\label{tab:perconvmix2}
\end{table}
\begin{table}[ht]
\centering
\begin{tabular}{llclc}
\hline
$N$ & $E_\infty[H]$ & rate for $E_\infty[H]$ & $E_\infty[U]$ & rate for $E_\infty[U]$\\
\hline
$10$ & $1.3308\times 10^{-3}$ & -- & $1.4064\times 10^{-3}$ & -- \\
$20$ & $2.6812\times 10^{-4}$ & $2.3113$ & $9.9809\times 10^{-5}$ & $3.8167$   \\
$40$ & $6.2091\times 10^{-5}$ & $2.1104$ & $9.1910\times 10^{-6}$ & $3.4409$   \\
$80$ & $1.5350\times 10^{-5}$ & $2.0161$ & $1.0541\times 10^{-6}$ & $3.1241$   \\
$160$ & $3.8219\times 10^{-6}$ & $2.0059$ & $1.2825\times 10^{-7}$ & $3.0390$  \\
$320$ & $9.5358\times 10^{-7}$ & $2.0029$ & $1.6047\times 10^{-8}$ & $2.9986$  \\
$640$ & $2.3815\times 10^{-7}$ &  $2.0014$ & $1.9751\times 10^{-9}$ & $3.0223$ \\
\hline 
\end{tabular}
\caption{Spatial errors and rates of convergence for the exact solution with $S_h=P^1$ and $S_u=P^2$ using the $L^\infty$ norm.}\label{tab:perconvmix3}
\end{table}

Finally,  Tables \ref{tab:perconvcub1}, \ref{tab:perconvcub2} and \ref{tab:perconvcub3} present the errors and convergence rates for the modified Galerkin method with cubic splines, i.e. with $S^3$ elements. In this case, the results are similar to those of the standard Galerkin method for the \acs{SGN} and Peregrine's system. Specifically, Table \ref{tab:perconvcub1} suggests a rate of convergence for $h$ close to $3.5$ which is in a good agreement with the estimate $\|h-\tilde{h}\|=O(\Delta x^{3.5}\sqrt{\ln (1/\Delta x)})$ suggested for the standard Galerkin method for  Peregrine's system \cite{AD2013}. Also Table \ref{tab:perconvcub1} suggests a rate of convergence for $h$ close to $4$, which agrees with the estimate $\|u-\tilde{u}\|=O(\Delta x^{4}\sqrt{\ln (1/\Delta x)})$ of \cite{AD2013}.

\begin{table}[ht]
\centering
\begin{tabular}{llclc}
\hline
$N$ & $E_0[H]$ & rate for $E_0[H]$ & $E_0[U]$ & rate for $E_0[U]$\\
\hline
$200$ & $0.5072\times 10^{-8}$ & -- &  $0.3101\times 10^{-10}$ &  -- \\
$250$ &  $0.2287\times 10^{-8}$ & $3.5697$ & $0.1270\times 10^{-10}$ & $3.9992$ \\
$300$ &  $0.1202\times 10^{-8}$ & $3.5256$ & $0.6127\times 10^{-11}$ & $3.9991$ \\
$350$ &  $0.6998\times 10^{-9}$ & $3.5125$ & $0.3307\times 10^{-11}$ & $3.9989$ \\
$400$ &  $0.4384\times 10^{-9}$ & $3.5019$ & $0.1939\times 10^{-11}$ & $3.9975$ \\
$450$ &  $0.2913\times 10^{-9}$ & $3.4710$ & $0.1210\times 10^{-11}$ & $4.0006$ \\
$500$ &  $0.2021 \times 10^{-9}$& $3.4679$ & $0.7968\times 10^{-12}$ & $3.9708$ \\
\hline 
\end{tabular}
\caption{Spatial errors and rates of convergence for the exact solution with $S^3$ finite elements using the $L^2$ norm.}\label{tab:perconvcub1}
\end{table}

\begin{table}[ht]
\centering
\begin{tabular}{llclc}
\hline
$N$ & $E_1[H]$ & rate for $E_1[H]$ & $E_1[U]$ & rate for $E_1[U]$\\
\hline
$200$ & $0.3406\times 10^{-5}$ & -- & $0.3920\times 10^{-7}$ & -- \\
$250$ & $0.1970\times 10^{-5}$ & $2.4531$ & $0.2008\times 10^{-7}$ & $2.9982$ \\
$300$ & $0.1265\times 10^{-5}$ & $2.4294$ & $0.1162\times 10^{-7}$ & $2.9986$ \\
$350$ & $0.8700\times 10^{-6}$ & $2.4299$ & $0.7322\times 10^{-8}$ & $2.9987$ \\
$400$ & $0.6286\times 10^{-6}$ & $2.4334$ & $0.4905\times 10^{-8}$ & $2.9988$ \\
$450$ & $0.4722\times 10^{-6}$ & $2.4291$ & $0.3446\times 10^{-8}$ & $2.9989$ \\
$500$ & $0.3654\times 19^{-6}$ & $2.4348$ & $0.2512\times 10^{-8}$ & $2.9989$ \\
\hline 
\end{tabular}
\caption{Spatial errors and rates of convergence for the exact solution with $S^3$ finite elements using the $H^1$ norm.}\label{tab:perconvcub2}
\end{table}

\begin{table}[ht]
\centering
\begin{tabular}{llclc}
\hline
$N$ & $E_2[H]$ & rate for $E_2[H]$ & $E_2[U]$ & rate for $E_2[U]$\\
\hline
$200$ & $0.2759\times 10^{-2}$ & -- & $0.5098\times 10^{-4}$ & -- \\
$250$ & $0.2033\times 10^{-2}$ & $1.3673$ & $0.3262\times 10^{-4}$ & $2.0010$ \\
$300$ & $0.1583\times 10^{-2}$ & $1.3715$ & $0.2265\times 10^{-4}$ & $2.0009$ \\
$350$ & $0.1278\times 10^{-2}$ & $1.3863$ & $0.1663\times 10^{-4}$ & $2.0006$ \\
$400$ & $0.1060\times 10^{-2}$ & $1.3999$ & $0.1273\times 10^{-4}$ & $2.0005$ \\
$450$ & $0.8987\times 10^{-3}$ & $1.4078$ & $0.1006\times 10^{-4}$ & $2.0004$ \\
$500$ & $0.7741\times 10^{-3}$ & $1.4168$ & $0.8151\times 10^{-5}$ & $2.0003$ \\
\hline 
\end{tabular}
\caption{Spatial errors and rates of convergence for the exact solution with $S^3$ finite elements using the $H^2$ norm.}\label{tab:perconvcub3}
\end{table}

\begin{table}[ht]
\centering
\begin{tabular}{llclc}
\hline
$N$ & $E_\infty[H]$ & rate for $E_\infty[H]$ & $E_\infty[U]$ & rate for $E_\infty[U]$\\
\hline
$200$ & $0.5211\times 10^{-7}$ & --  & $0.7553\times 10^{-10}$ & -- \\
$250$ & $0.2700\times 10^{-7}$ & $2.9460$ & $0.3111\times 10^{-10}$ & $3.9751$ \\
$300$ & $0.1577\times 10^{-7}$ & $2.9493$ & $0.1505\times 10^{-10}$ & $3.9812$ \\
$350$ & $0.1000\times 10^{-7}$ & $2.9559$ & $0.8145\times 10^{-11}$ & $3.9843$ \\
$400$ & $0.6732\times 10^{-8}$ & $2.9628$ & $0.4784\times 10^{-11}$ & $3.9857$ \\
$450$ & $0.4748\times 10^{-8}$ & $2.9649$ & $0.2990\times 10^{-11}$ & $3.9888$ \\
$500$ & $0.3471\times 10^{-8}$ & $2.9727$ & $0.1964\times 10^{-11}$ & $3.9870$ \\
\hline 
\end{tabular}
\caption{Spatial errors and rates of convergence for the exact solution with $S^3$ finite elements using the $L^\infty$ norm.}\label{tab:perconvcub4}
\end{table}
\begin{figure}[ht!]
  \centering
  \includegraphics[width=\columnwidth]{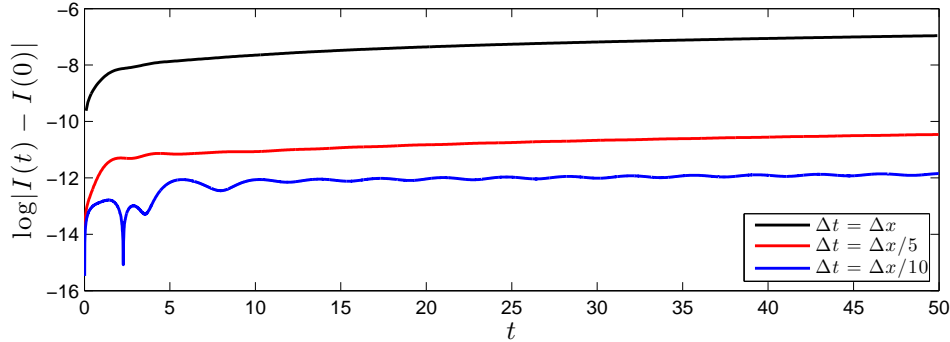}
  \caption{Sensitivity of energy conservation to the value of $\Delta t$ for fixed $\Delta x=0.1$}
  \label{fig:energy1}
\end{figure}
\begin{figure}[ht!]
  \centering
  \includegraphics[width=\columnwidth]{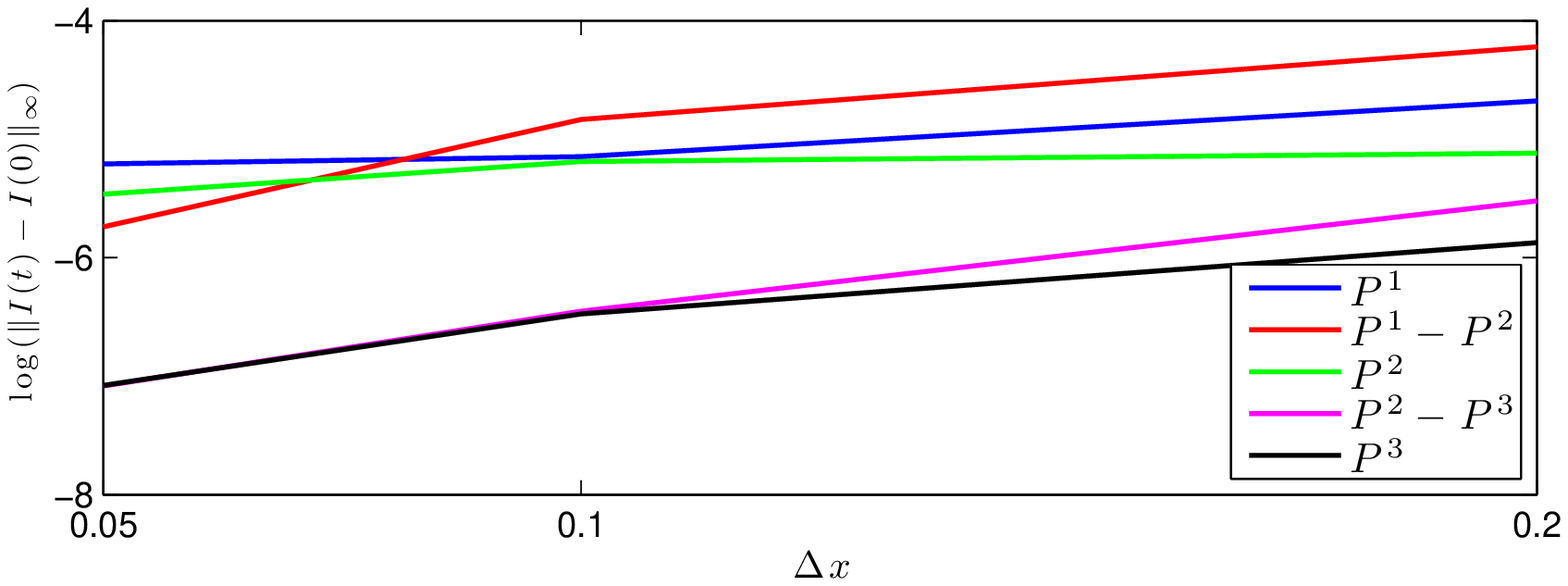}
  \caption{Sensitivity of energy conservation to the value of $\Delta x$ for small values of $\Delta t$}
  \label{fig:energy2}
\end{figure}

To study further the accuracy of the method we compute the conservation of the energy functional $I(t)$. As initial condition we use a solitary wave of amplitude $a=0.2$ with a smooth bottom given by the function $b(x)=1+0.1\sin(\pi x/2)$. Because the \acs{SGN} system together  with the wall boundary conditions does not conserve the energy functional, we consider the interval $[-100,100]$ and maximum time $T=50$ where the solution remained practically zero at the endpoints of the domain. We also used a fairly coarse grid with $\Delta x=0.1$. Due to the dissipative properties of the explicit Runge-Kutta method,  conservation of the invariant $I(t)$ is sensitive to the size of the timestep $\Delta t$. For example, when we used $\Delta x=\Delta t=0.1$ with elements in $S^3$, the invariant remained at the value $I(t)=0.314542$, when we used $\Delta t=0.01$, then the invariant remained at the value $I(t)= 0.31454249795$ conserving the digits shown. 

Although the energy functional $I(t)$ depends on the choice of $\Delta t$, the value of the Hamiltonian remains practically constant for any small value of $\Delta t<0.01$ and the error $|I(t)-I(0)|$ is of $O(10^{-12})$. Figure \ref{fig:energy1} shows the sensitivity of the invariant $I(t)$ to $\Delta t$ for fixed $\Delta x=0.1$ for $S^3$ finite elements.  Figure \ref{fig:energy2} shows the error of the energy functional as a function of $\Delta x$ for different finite element spaces. We observe that the convergence of the computed energy functional is connected to the convergence of the numerical method.  Although the explicit Runge-Kutta methods we used are not conservative (in the sense that they introduce  a small amount of numerical dissipation), and their contribution to the error in the conservation of energy is larger than the error coming from the spatial discretization of the same order,  it seems that the actual error is not important and upon choosing appropriate small values of $\Delta t$, it can be considered negligible. When we used $P^1$ and $P^2$ elements, then the error from the spatial discretisation is larger than the errors embedded by the temporal discretisation that are considered unimportant.

\section{Numerical experiments}\label{sec:numexp}

In this section we present a series of numerical experiments that serve as benchmarks to verify the accuracy of the modified Galerkin method, in simulations with variable bathymetry and wall boundary conditions. We tested all the numerical methods for all the numerical experiments and we report all the important differences in the numerical results. All the graphs contain the results obtained with $P^1$ elements. 

\subsection{Shoaling of solitary waves}

Shoaling of solitary waves and the nonlinear mechanisms behind the shoaling of a solitary wave have been studied theoretically and experimentally by various authors, \cite{GSSV1994,GSS1997,Syn1991,SynSkj1993}, and shoaling is closely related to the runup of solitary wave on a plain beach, \cite{Synolakis1987}. In \cite{SynSkj1993} after reviewing the derivation of  Green's law for the amplitude evolution of shoaling waves for \acs{BT} systems, that is,
\begin{equation}\label{eq:green}
\eta_{\max}/A\sim b(x)^{-1/4}\ , 
\end{equation}
 several experiments with shoaling of solitary waves are presented. In \cite{Syn1991}  Green's law \mbox{$\eta_{\max}/A\sim b(x)^{-1}$} for the Shallow Water Wave equations is derived. 
In \cite{SynSkj1993,Skj1987} four regions of shoaling were determined, where the general Green's law $\eta_{\max}\sim b(x)^{-\alpha}$ applies for different values of the parameter $\alpha$:~ The zone of gradual shoaling ($\alpha<1$), the zone of rapid shoaling $\alpha\ge 1$, the zone of rapid decay $\alpha<-1$ and the zone of gradual decay $\alpha\ge -1$.  Green's law for Boussinesq systems is valid in the zone of gradual and rapid shoaling, while Green's law for the shallow water system derived in \cite{Syn1991} applies in the zone of rapid shoaling. In this work, we consider first the shoaling of solitary waves on a plane beach of (mild) slope $1:35$. This experiment has been proposed by Grilli {\em et. al.} \cite{GSSV1994,GSS1997}. 

\begin{figure}[ht!]
  \centering
  \includegraphics[width=\columnwidth]{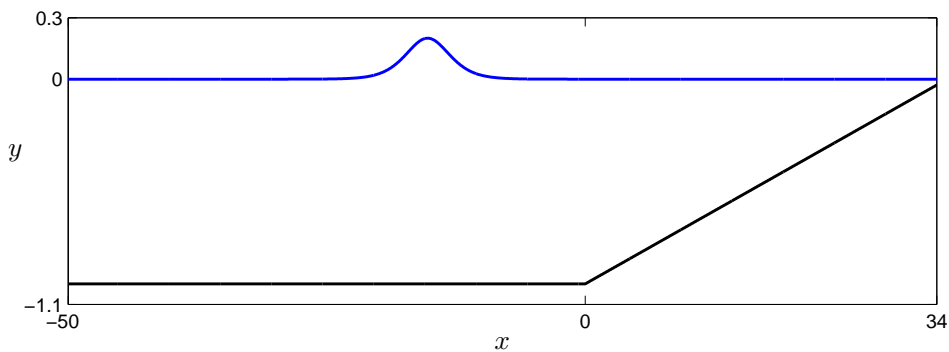}
  \caption{Sketch of the domain for the shoaling of solitary waves on a plain beach of slope $1:35$.}
  \label{fig:shoaling1}
\end{figure}

Here we use experimental data taken from \cite{GSSV1994} to compare with our numerical solution. For the numerical simulation we considered the domain shown in Figure \ref{fig:shoaling1}. For this experiment we take a uniform mesh on $[-100,34]$ with meshlength $\Delta x=0.1$. We also translate the solitary waves so that the crest amplitude is achieved at $x=-20.1171$. Because of the low regularity of the bottom at $x=0$, artifacts of the numerical method might appear especially in the case of $S^3$ elements. For this reason, we ensure the required regularity of the seafloor (required by the model equations), by approximating the bottom topography function with a quadratic polynomial near $x=0$. 

Although the \acs{SGN} system cannot model accurately the breaking of solitary waves, it appears that it models shoaling with higher accuracy than other Boussinesq models \cite{FBCR2015}, even very close to the breaking point. For example when $A=0.2$ the steep wave observed in the laboratory at the location called gauge 9 of \cite{GSSV1994} is approximated by a smooth solution in the \acs{SGN} system.  

\begin{figure}[ht!]
  \centering
  \includegraphics[width=\columnwidth]{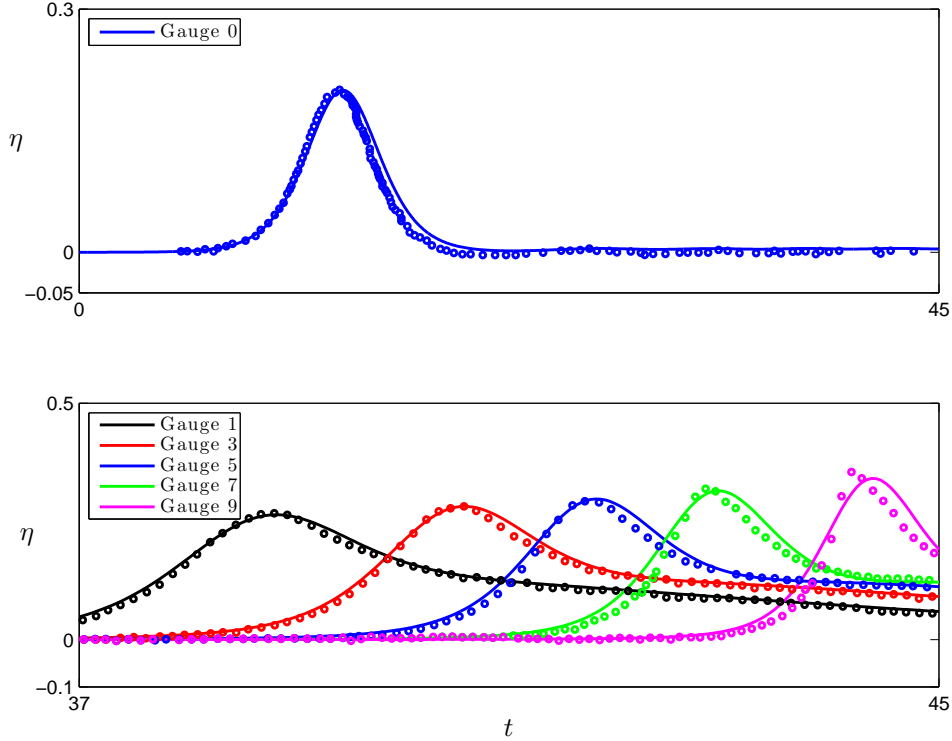}
  \caption{Solution on various wave gauges for the shoaling on a plane beach of slope $1:35$ of a solitary wave with $A=0.2$. Circles show experimental data, and lines show numerical solutions.}
  \label{fig:shoaling2}
\end{figure}

Next, we study the shoaling of solitary waves with normalized amplitude $A=0.1$, $0.15$, $0.2$ and $0.25$.
We monitor the numerical solution on the gauges (enumerated as in \cite{GSSV1994}) with number 0, 1, 3, 5, 7 and 9 located at $x=-5.0, 20.96, 22.55, 23.68, 24.68, 25.91$ respectively. The results with $P^1$ elements in the case of the shoaling of the solitary wave  with $A=0.2$ are shown in Figure \ref{fig:shoaling2}. The results with the other methods are almost identical and are not shown here.

\begin{figure}[ht!]
  \centering
  \includegraphics[width=\columnwidth]{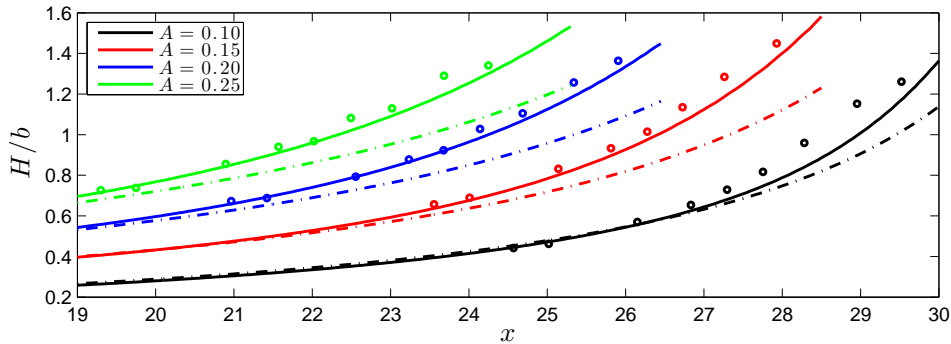}
  \caption{Comparison of the computed relative wave height $H/b$ with  experimental data of \cite{GSSV1994} for the shoaling of solitary waves on a plane beach of slope $1:35$, and with Green's law. '$-$': Numerical solutions, '$\circ$': Experimental data, '$-\cdot$': Green's law for Boussinesq systems}
  \label{fig:shoaling3}
\end{figure}
In these experiments, we also monitored the invariant $I(t)$ for values of $t$ small enough that the waves do not interact with the boundaries. In the case of  $S^3$ elements on a uniform grid with $\Delta x=0.1$, the values of the invariant up to time $t=30$ are given in the Table \ref{tab:shoalinv}. For the same experiments, the respective invariants were conserved to 4 decimal digits in the case of $P^1$ elements with $\Delta x=0.1$.
\begin{table}[ht]
\centering
\begin{tabular}{llclc}
\hline
$A$ & $I(t)$ \\
\hline
$0.10$ & $0.104058609813$\\
$0.15$ & $0.197139475070$\\
$0.20$ & $0.312548348249$\\
$0.25$ & $0.449208354485$ \\
\hline 
\end{tabular}
\caption{The conserved invariant $I(t)$ for up to $t=30$.}
\label{tab:shoalinv}
\end{table}
Finally, we computed the relative wave height defined as $H(x^\ast)/b(x^\ast)$ where $H(x^\ast)=\max_x\{|\zeta(x,t)|\}$ is the crest amplitude of the wave normalized by the local depth $b(x^\ast)$ evaluated at the same point $x^\ast$. Figure \ref{fig:shoaling3} presents the experimental and numerical data. We observe a very good agreement between the numerical and the experimental results. We note that the numerical results obtained by the numerical solution of the full water wave problem in \cite{GSSV1994} fit the experimental data (see e.g. Fig. 4 of \cite{GSSV1994}), as well as  the \acs{SGN} model. 

It is noted that the evolution of the maximum of the solution according to Green's law for Boussinesq systems, \cite{SynSkj1993}, represented in Figure \ref{fig:shoaling3} by a broken line is initially close to the numerical and experimental data, then underestimates the evolution of the maximum especially during the shoaling of large amplitude solitary waves. This indicates that the solitary waves are  in between the zone of rapid and gradual shoaling since they are nearly breaking waves. On the other hand  Green's law for the nonlinear shallow water system, \cite{Syn1991}, which is omitted here, overestimates the amplitude of the shoaling solitary wave, likely because of the simplifications inherent in the derivation of Green's law. These results verify the conclusions made in \cite{SynSkj1993} about the zone of gradual shoaling where a law that is similar to Green's law $\eta_{\max}\sim b^{-1/4}$ is valid. In our case we computed by experimentation that $\eta_{\max}\sim b^{-2/7}$ describes quite well the evolution of the amplitude of a shoaling solitary wave on a plane beach of slope $1:35$.

\subsection{Reflection of solitary waves on vertical wall}

We study also the reflection of the solitary waves on a vertical wall. In general, a vertical wall can be modeled by assuming that there is no flux through the vertical wall, i.e. the horizontal velocity of the fluid on the wall is $u=0$. The reflection of a solitary wave at a vertical wall is equivalent to the head-on collision of two counter-propagating solitary waves of the same shape. Because the head-on collision implies that at the center of the symmetric interaction $\eta_x=0$ then one might argue that the symmetric head-on collision is governed by different mathematical properties compared to the interactions of a solitary wave at a wall. 

In practice, during the interaction of a solitary wave with a wall the derivative $\eta_x$ on the boundary is negligible and so the additional condition seems to be satisfied and the reflection at the wall is apparently  equivalent to the head-on collision of two counter propagating waves, \cite{AD2012}. (Although the boundary condition $\eta_x=0$ is not necessary, its imposition does not require any modification to the numerical method or to the finite element spaces). The interaction of solitary waves with a vertical wall and the head-on collision of the solitary waves has been studied by theoretical and numerical means in \cite{Mirie1982}. In \cite{Mirie1982}, an asymptotic solution for the maximum runup has been derived, namely, if $\alpha$ is the normalized amplitude $A/b$ of the impinging wave, then the normalized maximum runup $\bar{R}_{\max}=R_{\max}/b$ is approximately 
\begin{equation}\label{eq:maxrup}
\bar{R}_{\max}\sim 2\alpha+\half \alpha^2+\half\alpha^3.
\end{equation}
In \cite{Mirie1982}, the \acs{SGN} system with horizontal bottom has been solved  numerically using a finite difference scheme. Finite difference schemes for the \acs{SGN} equations have certain disadvantages, including the introduction of numerical dispersion or dissipation. In addition it is necessary to use less physical boundary conditions on the wall, such as $h_x=u_{xx}=u-\frac{1}{3h}(h^3u_x)_x=0$, \cite{Mirie1982}. A fourth-order compact finite volume scheme was used to solve the \acs{SGN} with general bathymetry in \cite{CBB1, CBB2}. The derivation of the boundary conditions proposed in \cite{CBB2} follow analogous techniques with \cite{P1967} imposing also more boundary conditions for both $h$ and $u$ that approximate the wall boundary conditions. Other useful boundary conditions such as absorbing boundary conditions were constructed in \cite{CBB2} and can be implemented equally well with the present numerical model. 

The head-on collision of solitary waves in the full water wave problem has also been studied asymptotically, in \cite{SM1980}. Specifically, a similar formula to the one for the \acs{SGN} equations has been derived, and is given by the formula $2\alpha+\frac{1}{2} \alpha^2+\frac{3}{4}\alpha^3$. Other somewhat less accurate asymptotic approximations have been derived in \cite{PTGOP1999}. 
\begin{figure}[ht!]
  \centering
  \includegraphics[width=\columnwidth]{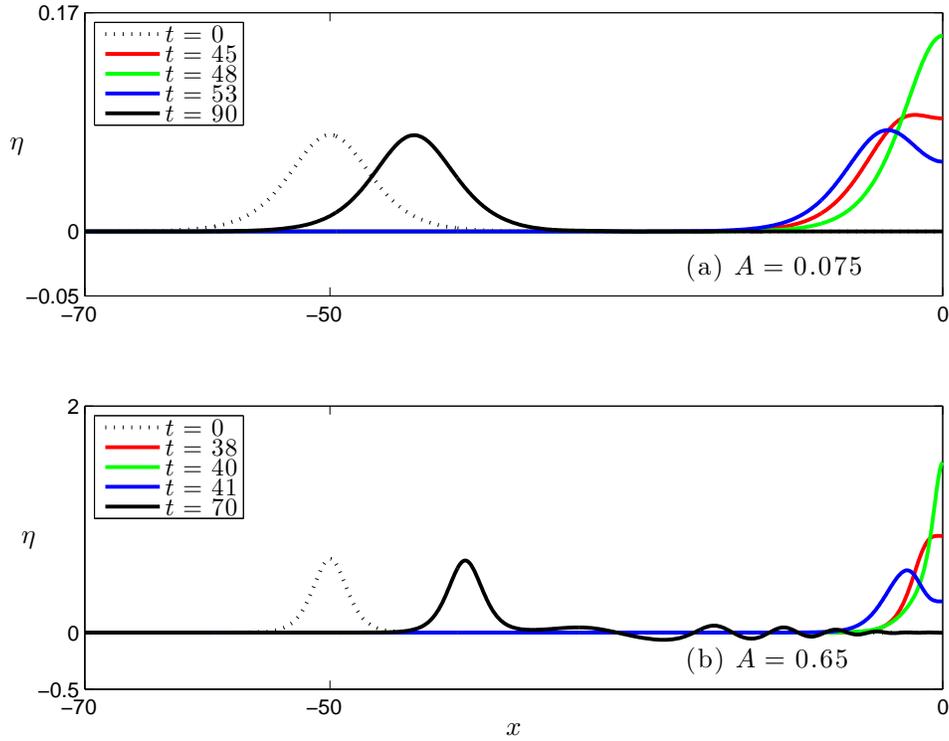}
  \caption{Reflection of solitary waves at a vertical wall located at $x=0$.}
  \label{fig:colis1}
\end{figure}
\begin{figure}[ht!]
  \centering
  \includegraphics[width=\columnwidth]{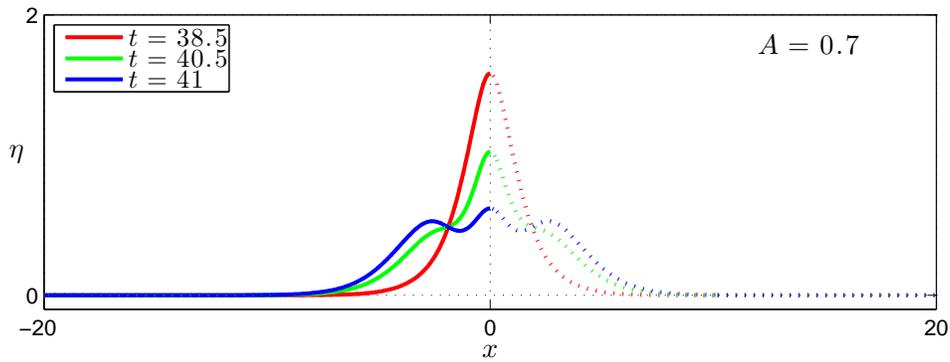}
  \caption{Reflection of solitary waves of normalized amplitude $A=0.7$ at a vertical wall located at $x=0$, plotted together with the head-on collision of two solitary waves of the same amplitude to show their equivalence.}
  \label{fig:colis2}
\end{figure}
Recent numerical experiments showed that \acs{BT} models do not describe accurately the reflection of large amplitude solitary waves at a vertical wall. Specifically,  in \cite{CWB1997} solving the full water wave equations using boundary element methods,  it was shown that large amplitude waves can achieve a higher maximum runup during a head-on collision with a vertical wall, than the predicted values of the respective asymptotic solutions of \cite{SM1980,Mirie1982} for solitary wave reflection at a vertical wall. These results have been verified in \cite{Chambarel2009,CKYHTC2015}.  Additionally, the formation of a residual jet was observed during the head-on collision of two large amplitude solitary waves, which was responsible for the large values of the maximum runup height.

To verify the ability of the numerical method to model accurately the reflection of solitary waves, we considered the case of a horizontal bottom $b(x)=-1$ for $x\in [-100,0]$ and solitary waves with amplitudes $A=0.075, 0.1, 0.15, \cdots, 0.7$. This case and these waves match the experiments which are presented in \cite{CWB1997} and which serve as benchmarks for the models in \cite{Chambarel2009,CKYHTC2015}.  We also consider a uniform grid on the interval $[-100,0]$ with $\Delta x=0.1$ while we translate the solitary waves so that their crest amplitude is at $x=-50$. Figure \ref{fig:colis1} shows the reflection of solitary waves with $A=0.075$ and $0.65$.
\begin{figure}[ht!]
  \centering
  \includegraphics[width=\columnwidth]{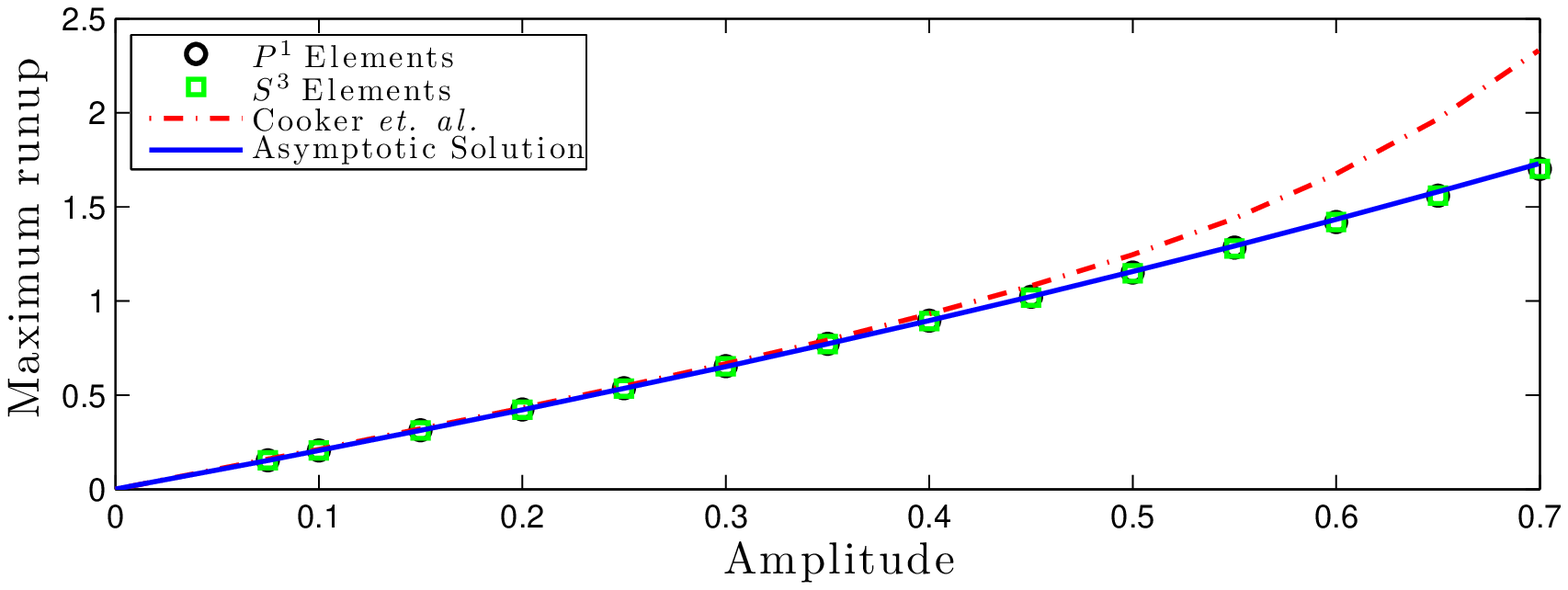}
  \caption{Maximum runup values during the reflection of solitary waves at a vertical wall located at $x=0$.}
  \label{fig:comrup}
\end{figure}

The reflection at a vertical wall of the extreme case of solitary waves of amplitude $A=0.7$ is presented in Figure \ref{fig:colis2}, together with  the solution of the head-on collision of two symmetric solitary waves. In this figure, we cannot observe any differences between the solutions of the reflected wave and of the colliding solitary waves (within  graphical accuracy). During this interaction a jet-like structure is visible, similar to those reported in \cite{Chambarel2009,CKYHTC2015,Maxworthy1976}. The jet formed in the case of the \acs{SGN} system is a solution of the mathematical model while the respective jet computed in \cite{Chambarel2009} is not. It is noted that no wave breaking is observed during these head-on collisions. Due to the inelastic interaction, the reflected solitary wave is followed by a dispersive tail. The magnitude of the dispersive tail depends on the amplitude of the colliding solitary wave. As it can be observed in Figure \ref{fig:colis1} the reflection of large solitary waves results in the generation of large dispersive tails. 

In \cite{Chambarel2009}, it was reported that the formation of a residual jet for the full water wave equations starts when the normalized amplitude of two incident solitary waves is larger than $0.60$. In the case of the \acs{SGN} equations, this jet formation was observed in reflections of solitary waves with amplitudes larger than $0.65$. A comparison of the numerical maximum runup values and the asymptotic formula (\ref{eq:maxrup}) in Figure \ref{fig:comrup} shows a very close match between the numerical solutions and the asymptotic results. In this figure, the numerical results of both the piecewise linear and the cubic spline finite elements are presented but no differences can be observed within graphical accuracy. Moreover, the numerical results of \cite{CWB1997} are compared with the numerical results obtained with the modified Galerkin method and we verify the difference in the head-on collision processes between the \acs{SGN} and the Euler equations, \cite{CWB1997,Chambarel2009}. This can be explained by noting that the jet cannot be modeled by the \acs{SGN} equations in the present form, since this artifact cannot be described by a smooth function but only by a parametric curve, \cite{Chambarel2009}.

\begin{remark}
In some cases it is considered that the reflection of a wave-train of more than one pulse can result in extreme runup values. Extreme wave runup on a vertical wall has been studied using periodic-boundary conditions and the head-on collision of wave-trains propagating in different directions in \cite{VCD2014,CDDD2013}. We repeated several of these experiments using half of the domain (required by the periodic code and with the wall boundary conditions described in this paper). The results were identical (except for the accuracy) to those obtained in \cite{VCD2014,CDDD2013} while our code  remained stable during the strong interactions of the waves with the wall. This verifies that the wall boundary conditions and the new numerical method can describe accurately strong interactions and extreme wave runup on a wall. As an indication of the results obtained we report here only the case of the reflection of a wave-train with three pulses ($N_w=3$, $\lambda_0=125 d$ in the notation of \cite{VCD2014}). In this case the normalized maximum runup observed was $R_{\max}/a_0\approx 5.545$.  
\end{remark}

\subsection{Reflection of shoaling waves}

\begin{figure}[ht!]
  \centering
  \includegraphics[width=\columnwidth]{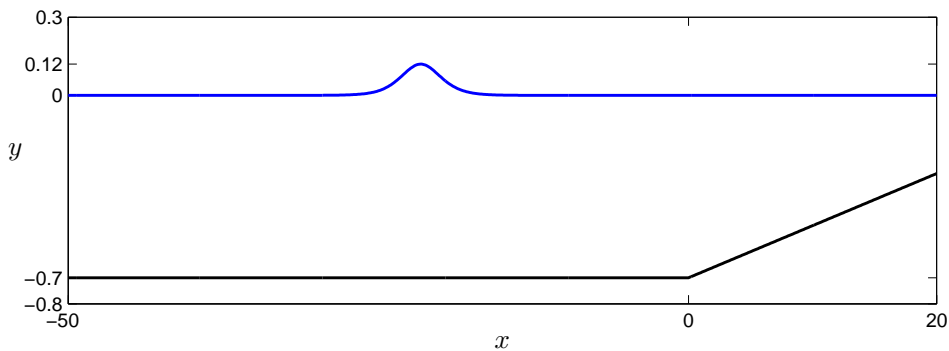}
  \caption{Sketch of the numerical experiment for the reflection at a vertical wall located at $x=20$ of a shoaling wave over a plain beach of slope $1:50$.}
  \label{fig:sketch2}
\end{figure}

In our final set of numerical experiments we consider the benchmark described in \cite{WB1999,Dodd1998}. Specifically, we consider the computational domain $[-100,20]$ (with $\Delta x=0.1$) and a bottom topography consisting of a horizontal seafloor in $[-100,0]$ and a plane beach of slope $1:50$ for $x\in[0,20]$. A sketch of the computational domain is presented in Figure \ref{fig:sketch2}. The vertical wall is located at $x=20$.

\begin{figure}[ht!]
  \centering
  \includegraphics[width=\columnwidth]{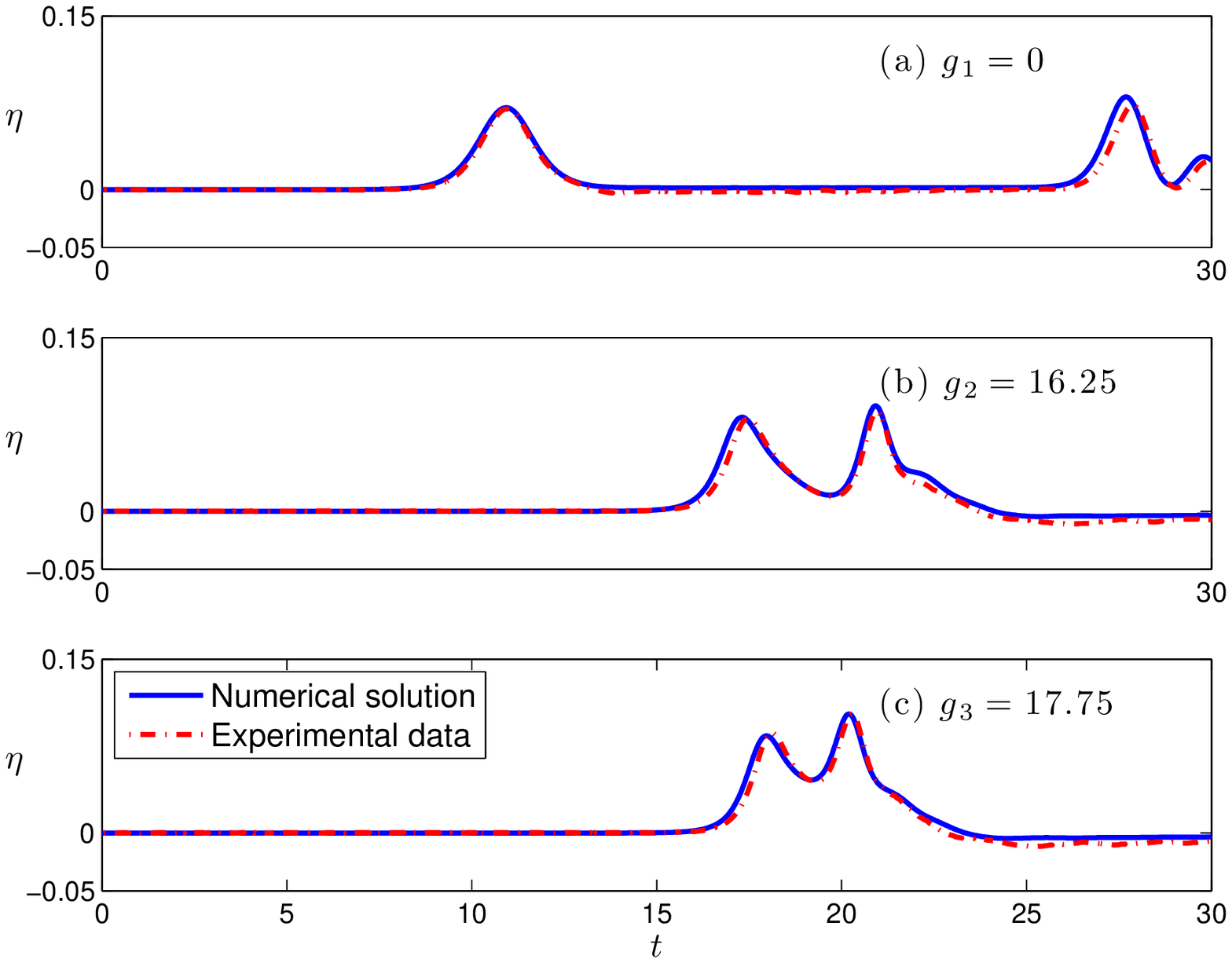}
  \caption{Reflection at a vertical wall located at $x=20$ of a shoaling wave over a plain beach of slope $1:50$. Initial solitary wave amplitude $A=0.07$.}
  \label{fig:figure6}
\end{figure}

\begin{figure}[ht!]
  \centering
  \includegraphics[width=\columnwidth]{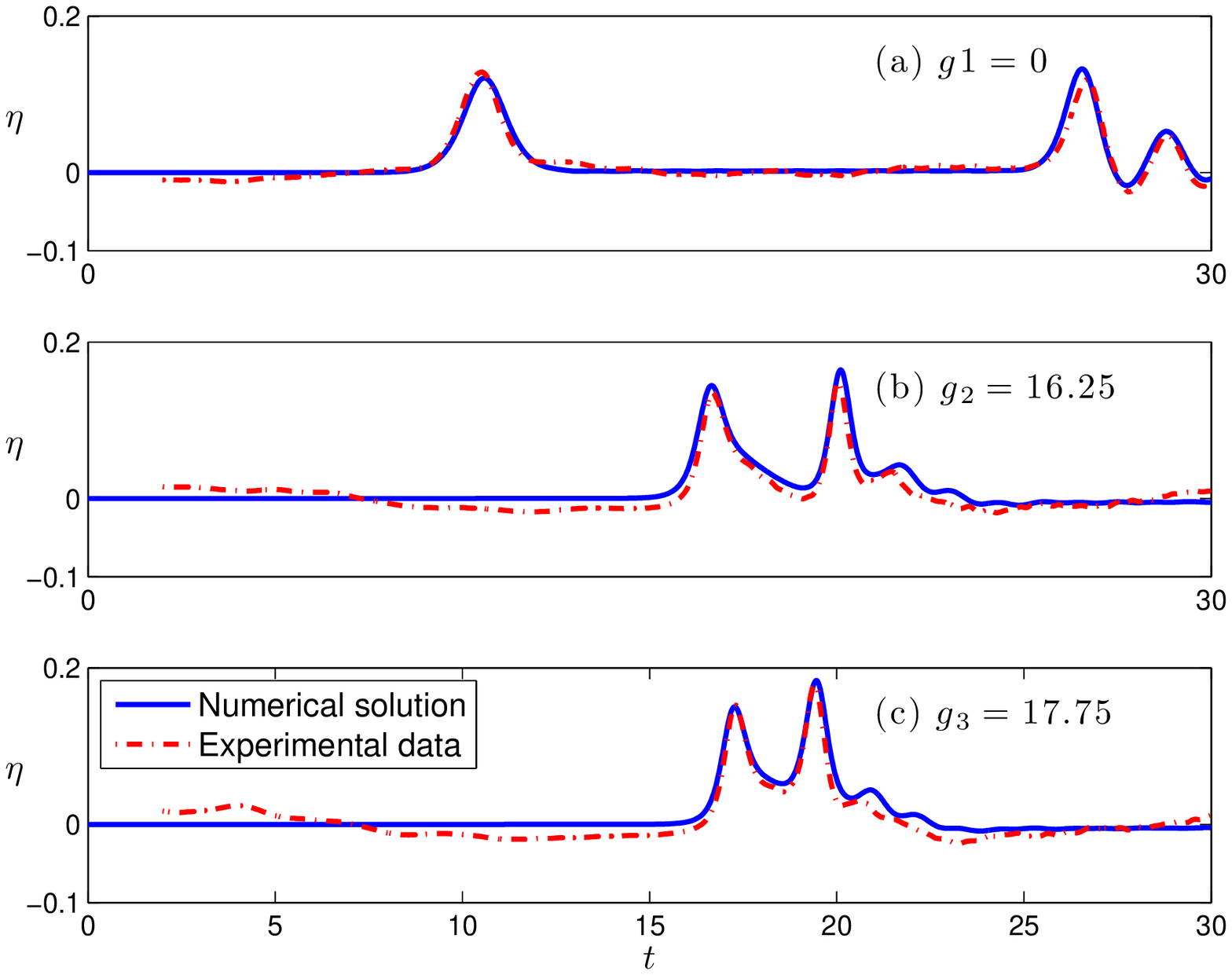}
  \caption{Reflection at a vertical wall located at $x=20$ of a shoaling wave  over a plain beach of slope $1:50$. Initial solitary wave amplitude $A=0.12$.}
  \label{fig:figure7}
\end{figure}

In this set of experiments we study the reflection of solitary waves at a vertical wall after the waves have climbed up the sloping beach. We consider two cases, one with a solitary wave of amplitude $A=0.07$ and another with $A=0.12$. The initial conditions have been translated so that the maximum of the crest is at $x=-30$. Figures \ref{fig:figure6} and \ref{fig:figure7} show experimental data recorded at three locations, by wave gauges $g_1$, $g_2$, and $g_3$, placed at $x=0$, $16.25$ and $17.75$. The comparison of numerical solutions with experimental data shows a better match from the present \acs{SGN} model than is obtained using any weakly nonlinear and weakly dispersive Boussinesq  system, \cite{FBCR2015,WB1999}. These numerical experiments verify the ability of the present numerical scheme to approximate with high accuracy and model fully nonlinear and weakly dispersive waves with wall boundary conditions at the endpoints of the computational domain.

In order to study the stability of the modified Galerkin method in more demanding situations, we consider the propagation of a solitary wave over a composite beach simulating the geometrical dimensions of Revere Beach, and its reflection by vertical wall. These experiments were conducted at the Coastal Engineering Laboratory of the U.S. Army Corps of Engineers, Vicksburg, Mississippi facility, \cite{Kanoglu1998} and serve as benchmarks for the reflection of nonbreaking, nearly breaking and breaking solitary waves by vertical wall. The composite beach consists of three piecewise linear segments while the bathymetry is constant  away from the beach and equal to $b_0=0.218~m$. The bathymetry can be realized by the function:
$$b(x)=\left\{\begin{array}{lr}
-0.218, & -11.77\leq x<15.04\\
1/53~x-0.5018, & 15.04 \leq x < 19.4\\
1/150~x-0.2650, & 19.4\leq x < 22.33\\
1/13~x-1.8340, &  22.33\leq x\leq 23.23
\end{array} \right. \ .$$

\begin{figure}[ht!]
  \centering
  \includegraphics[width=\columnwidth]{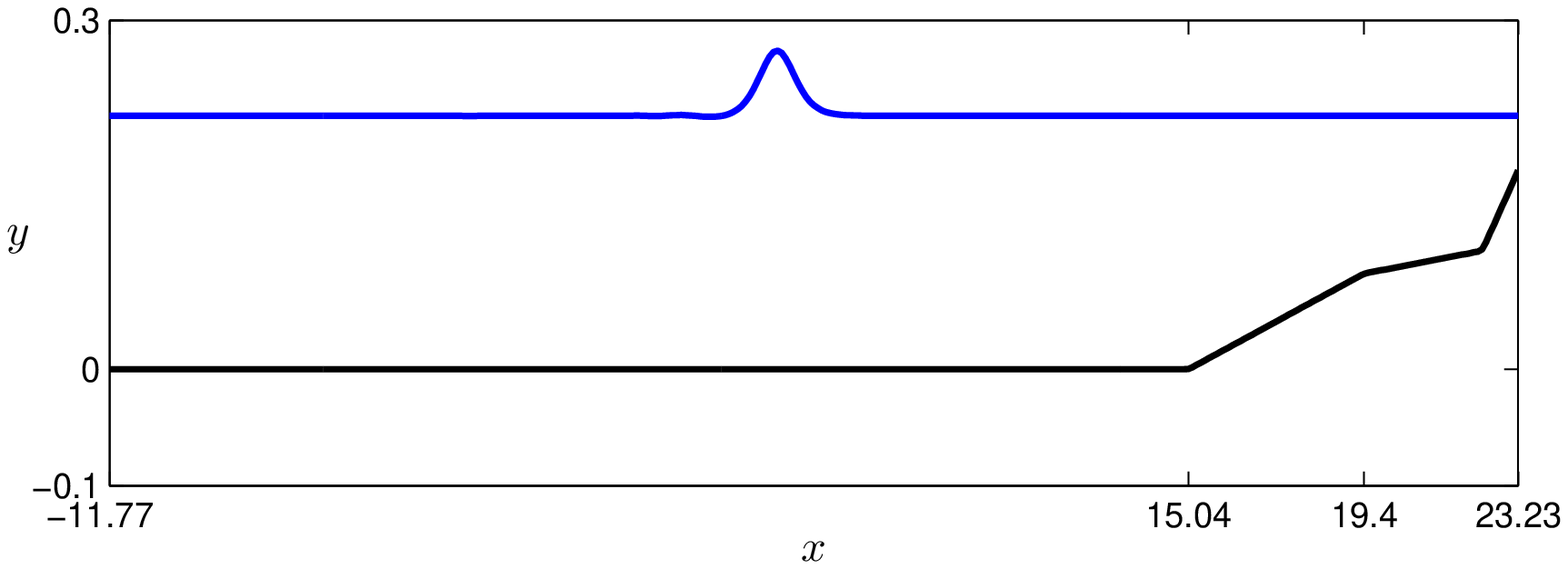}
  \caption{Sketch of the domain for the reflection of a solitary wave over a composite beach.}
  \label{fig:figure8b}
\end{figure}

The sketch of the domain is presented in Figure \ref{fig:figure8b}.
\begin{figure}[ht!]
  \centering
  \includegraphics[width=\columnwidth]{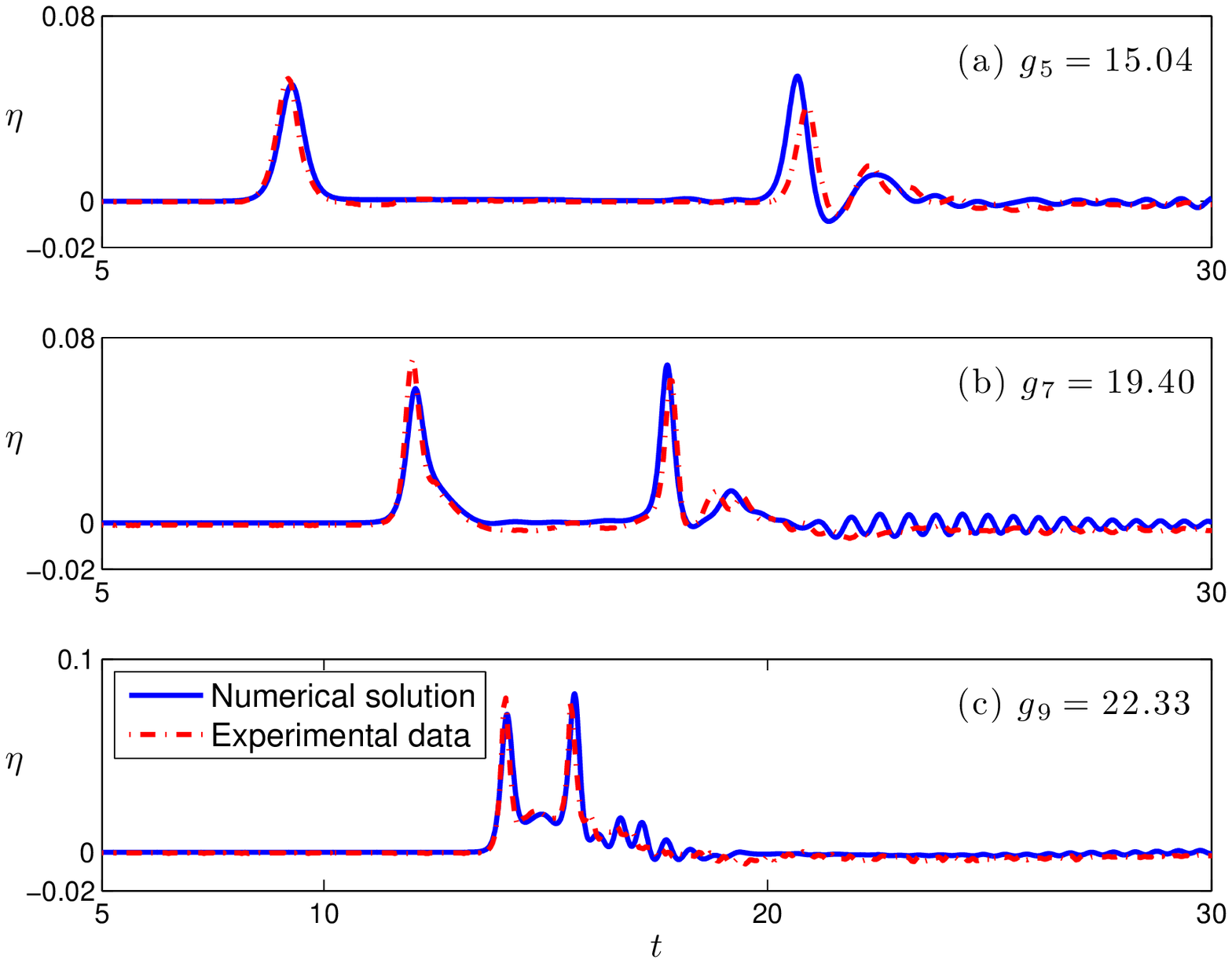}
  \caption{Reflection of a solitary wave over a composite beach. Vertical wall located at $23.23$. Initial solitary wave amplitude $A/b_0=0.3$.}
  \label{fig:figure9}
\end{figure}
Here, we considered three solitary waves with normalized amplitudes $A/b_0=0.05$, $0.3$ and $0.7$. We monitored the water depth  at $x$ locations  that correspond to gauges 5, 7 and 9 in  \cite{Kanoglu1998}. In this experiment, we used $\Delta x=0.1$ in $[-11.77,23.23]$. The normalized maximum runup computed at the location of the vertical wall in the first case was $R/b_0=0.122$ which is very close to the experimental value measured in the laboratory $R/b_0=0.13$. Moreover, the computed solution is very close to the experimental data at the gauges and we don't present these results here. It is noted that because a very small relative amplitude wave is involved in this case, its reflection can be modelled quite accurately, even by nondispersive models. The other two cases involve a nearly breaking and a breaking wave. Due to the steepness of the wave in the third case a wave breaking mechanism should be considered in order to approximate the solution in a stable manner. In the second experiment, the solitary wave is a nearly breaking wave. Although in this case the wave becomes very steep during shoaling, the numerical maximum runup computed was $0.46~m$, which is again very close to the experimentally recorded runup value $0.45~m$.  The solution at the wave gauges is presented in Figure \ref{fig:figure9}. The results can be improved by considering wave breaking mechanisms \cite{SSS2004,ACSS1996, Zelt1991} and Green-Naghdi models with improved dispersion characteristics such as those proposed in \cite{BN1996,DM2010,AC2013,AC2013b,LaMa2015}. It is noted that  in both cases only linear terms with second order derivatives are included locally around the regions where the waves become steep. These terms can be incorporated easily by any numerical method while the possible stiffness induced by the new terms  to the problem can be handled by taking smaller time-steps or by using a time integration method appropriate for stiff problems, \cite{ACSS1996}.

\section{Summary and Conclusions}\label{sec:summary}

We present a fully discrete numerical scheme for the \acf{SGN} system with wall boundary conditions. Semidiscretization of the model equations is based on a modification of the standard Galerkin / finite-element method that allows solutions in function spaces of low regularity. The time discretization is based on a fourth-order, four-stage, explicit Runge-Kutta method. A detailed computational study of the convergence properties of the numerical scheme shows that the method converges with similar convergence rates to those of  Peregrine's system. Some of the advantages of the new numerical method are:
\begin{itemize}
\item the high accuracy and the very good conservation properties;
\item the sparsity of the resulting linear systems;
\item the low complexity of the algorithm due to the use of low-order finite element spaces; and
\item its potential to be extended to the two-dimensional model equations.
\end{itemize}

In addition, we perform a series of very accurate numerical experiments to verify the efficacy of the numerical scheme in studies of shoaling and reflecting solitary waves. The numerical solutions are compared with available experimental data and theoretical approximations, whenever possible. The numerical model appears to be efficient and the match between numerical results, experimental data, and theoretical approximations is very satisfactory and shows better performance than other shallow water wave systems when there is no wave breaking. Wave breaking can be treated by following heuristic methodologies, such as discarding dispersive terms or adding new dissipative terms.

\section*{Acknowledgment}

The authors were supported by the Marsden Fund administered by the Royal Society of New Zealand.

\bibliographystyle{plain}
\bibliography{biblio}

\end{document}